\documentclass[11pt]{article}
\usepackage{smile}

\usepackage{fullpage}
\usepackage{lscape}
\usepackage{bigints}
\usepackage{upgreek}
\usepackage{framed}
\usepackage{mdframed}
\usepackage{enumerate}
\usepackage[inline]{enumitem}
\usepackage[T1]{fontenc}
\usepackage{moresize}
\usepackage{bm}
\usepackage{bbm}
\usepackage{dsfont}
\usepackage{amsmath}
\usepackage{amssymb}
\usepackage{amsthm}
\usepackage{amsfonts}
\usepackage{stmaryrd}
\usepackage{array}
\usepackage{mathrsfs}
\usepackage{mathtools} 
\usepackage{extarrows}
\usepackage{stackrel}
\usepackage{relsize,exscale}
\usepackage{scalerel}
\usepackage[nodisplayskipstretch]{setspace}
\usepackage{color}
\usepackage[usenames,dvipsnames]{xcolor}
\usepackage{cancel}
\usepackage{soul}
\usepackage{xfrac}
\usepackage{siunitx}
\usepackage{graphicx}
\usepackage{float}
\usepackage{rotating}
\usepackage{subcaption}
\usepackage{overpic}
\usepackage[all]{xy}
\usepackage{tikz}
\usetikzlibrary{arrows,matrix,positioning,calc,automata,patterns}
\usepackage{booktabs}
\usepackage{dcolumn}
\usepackage{multirow}
\usepackage{diagbox}
\usepackage{tabularx}
\usepackage{verbatim}
\usepackage{listings}
\usepackage[ruled,vlined]{algorithm2e}
\usepackage{fancyvrb}
\usepackage{hyperref}
\usepackage[round]{natbib}
\usepackage{sectsty}
\usepackage{etoolbox}
\usepackage[capitalize,noabbrev]{cleveref}
\crefname{assumption}{Assumption}{Assumptions}
\newcounter{step}
\AtBeginEnvironment{proof}{\setcounter{step}{0}}
\crefname{step}{Step}{Steps}

\hypersetup{
    bookmarks=true,         
    unicode=false,          
    pdftoolbar=true,        
    pdfmenubar=true,        
    pdffitwindow=false,     
    pdfstartview={FitH},    
    pdftitle={My title},    
    pdfauthor={Author},     
    pdfsubject={Subject},   
    pdfcreator={Creator},   
    pdfproducer={Producer}, 
    pdfkeywords={key1, key2}, 
    pdfnewwindow=true,      
    colorlinks=true,        
    linkcolor=blue,         
    citecolor=blue,         
    filecolor=blue,         
    urlcolor=cyan           
}

\usepackage{stackengine}
\stackMath
\newcommand\tenq[2][1]{%
\def\useanchorwidth{T}%
\ifnum#1>1%
\stackunder[0pt]{\tenq[\numexpr#1-1\relax]{#2}}{\!\scriptscriptstyle\thicksim}%
\else%
\stackunder[1pt]{#2}{\!\scriptstyle\thicksim}%
\fi%
}

\makeatletter
\DeclareRobustCommand\widecheck[1]{{\mathpalette\@widecheck{#1}}}
\def\@widecheck#1#2{%
    \setbox\z@\hbox{\m@th$#1#2$}%
    \setbox\tw@\hbox{\m@th$#1%
       \widehat{%
          \vrule\@width\z@\@height\ht\z@
          \vrule\@height\z@\@width\wd\z@}$}%
    \dp\tw@-\ht\z@
    \@tempdima\ht\z@ \advance\@tempdima2\ht\tw@ \divide\@tempdima\thr@@
    \setbox\tw@\hbox{%
       \raise\@tempdima\hbox{\scalebox{1}[-1]{\lower\@tempdima\box
\tw@}}}%
    {\ooalign{\box\tw@ \cr \box\z@}}}
\makeatother

\newcommand{\var}{\mathrm{Var}}
\newcommand{\Diag}{\mathrm{Diag}}

\newcommand{\Step}[2][]{%
  \refstepcounter{step}%
  \paragraph*{Step \thestep: #2.}%
  \if\relax\detokenize{#1}\relax\else\label{#1}\fi
}

\numberwithin{equation}{section}
\theoremstyle{plain}

\newtheorem{proposition}{Proposition}[section]
\newtheorem{assumption}{Assumption}[section]

\newtheorem{thm}{Theorem}[section]  
\newtheorem{cor}{Corollary}[section]
\newtheorem{lem}{Lemma}[section]

\newtheorem{innercustomgeneric}{\customgenericname}
\providecommand{\customgenericname}{}
\newcommand{\newcustomtheorem}[2]{%
  \newenvironment{#1}[1]
  {%
   \renewcommand\customgenericname{#2}%
   \renewcommand\theinnercustomgeneric{##1}%
   \innercustomgeneric
  }
  {\endinnercustomgeneric}
}
\newcustomtheorem{customdefinition}{Definition}
\newcustomtheorem{customdefinitions}{Definitions}
\newcustomtheorem{customtheorem}{Theorem}
\newcustomtheorem{customassumption}{Assumption}
\newcustomtheorem{customlemma}{Lemma}
\newcustomtheorem{customexample}{Example}
\theoremstyle{definition}
\newtheorem{definition}{Definition}[section]

\newtheorem{remark}{Remark}[section]

\usepackage{enumitem}
\makeatletter
\newcommand{\mylabel}[2]{#2\def\@currentlabel{#2}\label{#1}}
\makeatother

\graphicspath{{./fig3/}}

\newcommand{\zd}[1]{}
\newcommand{\nb}[1]{}

\newcommand{\R}{\mathbb{R}}  

\newcommand{\N}{\mathbb{N}}

\newcommand{\E}{\mathbb{E}}
\renewcommand{\P}{\mathbb{P}}

\renewcommand\ml{\mathcal}
\newcommand\tp{\intercal}
\newcommand\wt{\widetilde}

\newcommand{\rw}{\rightarrow}

\renewcommand{\l}{\left}
\renewcommand{\r}{\right}
\renewcommand{\phi}{\varphi}
\renewcommand{\epsilon}{\varepsilon}

\newcommand{\be}{\begin{equation}}
\newcommand{\ee}{\end{equation}}
\newcommand{\bea}{\begin{equation}\begin{aligned}}
\newcommand{\eea}{\end{aligned}\end{equation}}
\renewcommand{\tr}{\operatorname{tr}}
\renewcommand{\var}{\operatorname{Var}}
\renewcommand\ind{\text{\large 1\hskip-0.29em I}}

\allowdisplaybreaks

\begin{document}

\setlength{\abovedisplayskip}{5pt}
\setlength{\belowdisplayskip}{5pt}
\setlength{\abovedisplayshortskip}{5pt}
\setlength{\belowdisplayshortskip}{5pt}
\hypersetup{colorlinks,breaklinks,urlcolor=blue,linkcolor=blue}

\title{\LARGE Limiting spectral distributions of large consistent rank correlation matrices}

\author{Zhaorui Dong\thanks{School of Data Science, The Chinese University of Hong Kong, Shenzhen, China; {\tt zhaoruidong@link.cuhk.edu.cn}},~~~
Fang Han\thanks{Department of Statistics, University of Washington, Seattle, WA 98195, USA; {\tt fanghan@uw.edu}},
	     ~~and~~
Jianfeng Yao\thanks{School of Data Science, The Chinese University of Hong Kong, Shenzhen, China; {\tt jeffyao@cuhk.edu.cn}}	    
}

\date{}

\maketitle

\vspace{-1em}

\begin{abstract}
We study random matrices whose entries are obtained by applying consistent rank correlations, such as Hoeffding's $D$, pairwise to a high-dimensional random vector with mutually independent components. Prior work has shown that, in the proportional high-dimensional regime, the empirical spectral distributions of large Kendall's tau and Spearman's rho matrices converge weakly almost surely to the Marchenko--Pastur law. By contrast, we prove that for consistent rank correlations such as Hoeffding's $D$, the limiting spectral distribution is given by the semicircle law. Our result thus generalizes a recent work of Dong, Han, and Yao (2025), who considered the special case of Chatterjee's rank correlation and established the first semicircle law for a large correlation matrix in the proportional regime.
\end{abstract}


\section{Introduction}

For each $n=1,2,\ldots$, let $\mX_1,\ldots,\mX_n$ be independent and identically distributed $p$-dimensional random vectors, where
\[
\mX_i=(X_{i,1},\ldots,X_{i,p})^\tp,\qquad i\in[n]:=\{1,2,\ldots,n\}.
\]
We work under the assumption that
\begin{align}\label{eq:assumption}
X_{1,1},\ldots,X_{1,p}\text{ are mutually independent and continuous}.
\end{align}
Our interest lies in the proportional high-dimensional regime, in which the distribution of $\mX_1$ may depend on $n$ and
\begin{align}\label{eq:asymptotics}
p=p(n),\qquad n\rw\infty,\qquad \frac{p}{n}\to\gamma\in(0,\infty).
\end{align}

\subsection{Consistent rank correlations}\label{sec:crc}

Our focus is on correlation matrices whose entries are rank-based U-statistics capable of detecting nonlinear and non-monotone dependence between pairs $(X_{1,j},X_{1,k})$. More specifically, we study $p\times p$ matrices 
\[
\widehat{\fR}_n=[\widehat{R}_{jk}]_{j,k\in[p]}, 
\]
which we call {\it consistent rank correlation matrices}, with entries given by
\begin{align}\label{eq:cons-rank-corr}
\widehat{R}_{jk}=
\begin{cases}
1, & \text{if } j=k,\\[0.3em]
\binom{n}{m}^{-1}
\sum_{1\le i_1<\cdots<i_m\le n}
h\Big(
\binom{X_{i_1,j}}{X_{i_1,k}},
\ldots,
\binom{X_{i_m,j}}{X_{i_m,k}}
\Big), & \text{if } j\ne k.
\end{cases}
\end{align}
Here $h:(\R^2)^m\rw\R$ is a fixed and {\it symmetric} kernel of order $m$ satisfying the following two critical properties:
\begin{enumerate}[label=(\roman*)]
\item {\bf Rank-based.} The kernel $h$ depends on the inputs only through the rankings of the data; see Definition \ref{def:rank-based} for a precise definition.

\item {\bf Consistent.} If $\mZ=(Z_1,Z_2)$ is Lebesgue absolutely continuous and $\mZ_1,\ldots,\mZ_m$ are $m$ independent copies of $\mZ$, then
\begin{align*}
\E\big[h(\mZ_1,\ldots,\mZ_m)\big]=0
\quad\text{if and only if}\quad
Z_1 \text{ is independent of } Z_2.
\end{align*}
\end{enumerate}

We call statistics of the form \eqref{eq:cons-rank-corr} {\it consistent rank correlations}, since they give rise to rank-based tests of independence that are consistent against general alternatives \citep{shi2020power}. In particular, \cite{MR4185806} established a number of results on their statistical properties, including the limiting distribution of
\[
\max_{j\ne k\in[p]} |\widehat{R}_{jk}|
\]
under \eqref{eq:assumption} and \eqref{eq:asymptotics}. See also \cite{shi2020power} for the power of independence tests using these consistent rank correlations, as well as \cite{MR3737306} for some related results on Spearman and Kendall correlation matrices.

Prominent examples of consistent rank correlations of the form \eqref{eq:cons-rank-corr} include Hoeffding's $D$ \citep{MR0029139}, Blum--Kiefer--Rosenblatt's $R$ \citep{MR0125690}, and Bergsma--Dassios--Yanagimoto's $\tau^*$ \citep{MR3178526,yanagimoto1970measures}. Specifically, for arbitrary $\mz_j=(z_{j,1},z_{j,2})^\tp\in\mathbb{R}^2$, $j=1,2,\cdots$, these statistics correspond to \eqref{eq:cons-rank-corr} with kernels given by
\begin{align}
&\text{\bf(Hoeffding's $D$)}\qquad h_D(\mz_1,\ldots,\mz_5) \notag\\
&:= \frac{1}{16}\sum_{(i_1,\ldots,i_5)\in\cP_5}
\Big[\big\{\ind(z_{i_1,1}\le z_{i_5,1})-\ind(z_{i_2,1}\le z_{i_5,1})\big\}
\big\{\ind(z_{i_3,1}\le z_{i_5,1})-\ind(z_{i_4,1}\le z_{i_5,1})\big\}\Big] \notag \\
&\qquad \times
\Big[\big\{\ind(z_{i_1,2}\le z_{i_5,2})-\ind(z_{i_2,2}\le z_{i_5,2})\big\}
\big\{\ind(z_{i_3,2}\le z_{i_5,2})-\ind(z_{i_4,2}\le z_{i_5,2})\big\}\Big], \label{eq:hoeffding} \\[0.25cm]
&\text{\bf(Blum--Kiefer--Rosenblatt's $R$)}\qquad h_R(\mz_1,\ldots,\mz_6) \notag\\\
&:= \frac{1}{32}\sum_{(i_1,\ldots,i_6)\in\cP_6}
\big[\ind(z_{i_1,1}\le z_{i_5,1})-\ind(z_{i_2,1}\le z_{i_5,1})\big]
\big[\ind(z_{i_3,1}\le z_{i_5,1})-\ind(z_{i_4,1}\le z_{i_5,1})\big] \notag\\\
&\qquad \times
\big[\ind(z_{i_1,2}\le z_{i_6,2})-\ind(z_{i_2,2}\le z_{i_6,2})\big]
\big[\ind(z_{i_3,2}\le z_{i_6,2})-\ind(z_{i_4,2}\le z_{i_6,2})\big], \label{eq:bkr}\\
&\text{and} \notag\\\
&\text{\bf(Bergsma--Dassios--Yanagimoto's $\tau^*$)}\qquad h_{\tau^*}(\mz_1,\ldots,\mz_4) \notag\\\
&:= \frac{1}{16}\sum_{(i_1,\ldots,i_4)\in\cP_4}
\Big\{
\ind(z_{i_1,1},z_{i_3,1}<z_{i_2,1},z_{i_4,1})
+\ind(z_{i_2,1},z_{i_4,1}<z_{i_1,1},z_{i_3,1})\notag \\
&\qquad\qquad
-\ind(z_{i_1,1},z_{i_4,1}<z_{i_2,1},z_{i_3,1})
-\ind(z_{i_2,1},z_{i_3,1}<z_{i_1,1},z_{i_4,1})
\Big\} \notag\\\
&\qquad \times
\Big\{
\ind(z_{i_1,2},z_{i_3,2}<z_{i_2,2},z_{i_4,2})
+\ind(z_{i_2,2},z_{i_4,2}<z_{i_1,2},z_{i_3,2})\notag \\
&\qquad\qquad
-\ind(z_{i_1,2},z_{i_4,2}<z_{i_2,2},z_{i_3,2})
-\ind(z_{i_2,2},z_{i_3,2}<z_{i_1,2},z_{i_4,2})
\Big\}. \label{eq:bdy}
\end{align}
Here, for each positive integer $m$, $\cP_m$ denotes the permutation group on $[m]$ and $\ind(\cdot)$ the indicator function.

In contrast to \cite{MR4185806}, the objective of the present paper is to study, under \eqref{eq:assumption} and \eqref{eq:asymptotics}, the spectrum of the standardized correlation matrix
\[
\widehat{\fW}_n=\sqrt{n}\,(\widehat{\fR}_n-\fI_p),
\]
where $\fI_p$ denotes the $p\times p$ identity matrix.

\subsection{Motivation}

Our work contributes to a growing body of literature on the spectral analysis of large random matrices arising in statistics and machine learning, with particular emphasis on those whose limiting spectral distributions (LSDs) deviate in an essential way from the Marchenko-Pastur (MP)  law.

To place our results in context, we briefly review the literature on LSD theory for large random matrices. Under the proportional regime \eqref{eq:asymptotics}, the seminal works of \cite{wigner1958distribution} and \cite{marvcenko1967distribution} established the two fundamental limiting laws in random matrix theory, namely, the semicircle law and the MP law. More precisely, the empirical spectral distributions (ESDs) of Wigner matrices and sample covariance matrices converge to these two limits, respectively. It was later discovered that the LSDs of a variety of correlation matrices of statistical interest, including large Pearson, Spearman, and Kendall correlation matrices (the latest up to an appropriate rescaling), are likewise governed by the MP law \citep{jiang2004limiting,bai2008large,bandeira2017marvcenko}. These results were further sharpened by \citet{dong:yao:2025} and \citet{doernemann:heiny:2025}, who derived necessary and sufficient conditions for the MP law to hold for Pearson correlation matrices.

These developments naturally raise the following question: under \eqref{eq:asymptotics}, are there statistically meaningful classes of large random matrices whose LSDs are not given by the MP law? A partial answer was provided by \cite{jh:yao:heavytailed:aos}, who studied large Pearson correlation matrices in the heavy-tailed setting and showed that their LSD is given by a convolution of the MP law with a heavy-tailed component. More recently, \cite{dong2025spectralanalysislargedimensional} discovered that the LSD of a Chatterjee rank correlation matrix, whose entries are given by pairwise Chatterjee rank correlations \citep{chatterjee2020new,lin2022limit}, is in fact the semicircle law. 

It is also worth noting that a recent and active line of work in random matrix theory investigates limiting spectral behavior beyond the proportional regime \eqref{eq:asymptotics}. In particular, when $p/n\to 0$, \cite{BaiYin1988SC} and \cite{dornemann2025ties} showed that the LSDs of Pearson's covariance matrix and the correlation matrices of Spearman and Kendall (the latest again up to suitable rescaling), much like the phenomenon discovered in \cite{dong2025spectralanalysislargedimensional}, are in fact governed by the semicircle law. In a related direction, \cite{fan2019spectral}, \cite{lu2025equivalence}, and \cite{dubova2023universality}, among many others, studied the LSDs of random inner-product kernel matrices, which may be viewed as generalized covariance matrices, under polynomial scaling regimes. Their results show that the limiting law in this case is given by a free additive convolution of the MP law and the semicircle law.

The present paper contributes to this growing literature on non-MP limits. Our main motivation comes from \cite{dong2025spectralanalysislargedimensional}, which, to the best of our knowledge, provides by far the {\it only} non-MP limit theorem for a large correlation matrix in the proportional regime \eqref{eq:asymptotics}. More specifically, we extend the result therein from Chatterjee's rank correlation to a broader class of consistent rank correlation coefficients, and in doing so provide a unified perspective that also interestingly connects with the works of \cite{BaiYin1988SC}, \cite{dornemann2025ties}, \cite{fan2019spectral}, \cite{lu2025equivalence}, and \cite{dubova2023universality}.


\subsection{Our results}

Consider a $p\times p$ (possibly random) matrix $\fA=\fA_p$ with real eigenvalues
\[
\lambda_1(\fA)\ge \lambda_2(\fA)\ge \cdots \ge \lambda_p(\fA).
\]
The ESD of $\fA$ is defined as the normalized counting measure of its eigenvalues
\[
F^{\fA}:=\frac{1}{p}\sum_{i=1}^p \delta_{\lambda_i(\fA)},
\]
where $\delta_x$ denotes the Dirac measure at $x$.

Next, let $\operatorname{W}(r)$ denote Wigner's semicircle law with center $0$ and radius $r>0$, namely, the probability distribution with Lebesgue density
\[
\rho(x;r)=\frac{2}{\pi r^2}\sqrt{(r^2-x^2)_+},
\]
where $(x)_+:=x\ind(x\ge 0)$ denotes the positive part of $x$.

Lastly, denote 
\[
\widehat{\fR}_n^D, ~~\widehat{\fR}_n^R,~~ {\rm and}~~ \widehat{\fR}_n^{\tau^*} 
\]
to be the sample correlation matrices based on Hoeffding's $D$, Blum--Kiefer--Rosenblatt's $R$, and Bergsma--Dassios--Yanagimoto's $\tau^*$, respectively, obtained by choosing the kernel $h(\cdot)$ in Section~\ref{sec:crc} to be $h_D$, $h_R$, and $h_{\tau^*}$ introduced in \eqref{eq:hoeffding}, \eqref{eq:bkr}, and \eqref{eq:bdy}, respectively. 

Our main theorem, presented in Section~\ref{sec:main}, yields the following {\it distribution-free} semicircle laws for these matrices.

\begin{cor}\label{cor:crc}
Under \eqref{eq:assumption} and \eqref{eq:asymptotics}, we have almost surely,
\begin{align*}
&F^{\sqrt{n}\Big(\widehat{\fR}_n^D-\fI_p\Big)} \Rightarrow \operatorname{W}(2\sqrt{2\gamma}/3),\\
&F^{\sqrt{n}\Big(\widehat{\fR}_n^R-\fI_p\Big)} \Rightarrow \operatorname{W}(2\sqrt{2\gamma}),\\
\text{and}\qquad
&F^{\sqrt{n}\Big(\widehat{\fR}_n^{\tau^*}-\fI_p\Big)} \Rightarrow \operatorname{W}(6\sqrt{2\gamma}/5).
\end{align*}
Here $\Rightarrow$ denotes weak convergence.
\end{cor}

\begin{figure}[t]
  \centering
  \begin{subfigure}[t]{0.32\textwidth}
    \centering
    \includegraphics[width=\linewidth]{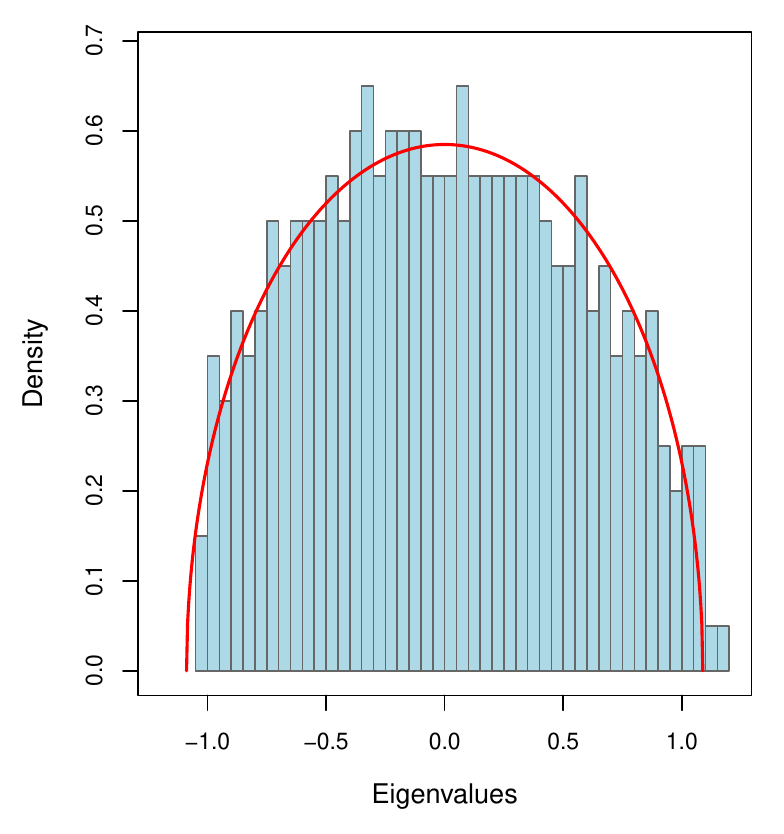}
    \caption{ESD of $\sqrt{n}\Big(\widehat{\fR}_n^D-\fI_p\Big)$.}
    \label{fig:sub-a}
  \end{subfigure}\hfill
  \begin{subfigure}[t]{0.32\textwidth}
    \centering
    \includegraphics[width=\linewidth]{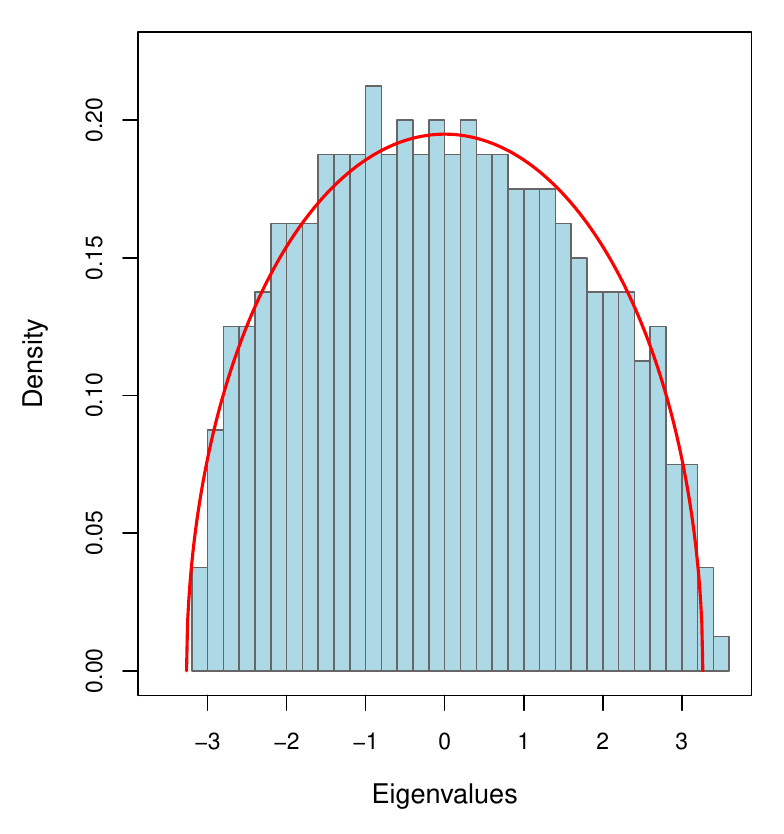}
    \caption{ESD of $\sqrt{n}\Big(\widehat{\fR}_n^{R}-\fI_p\Big)$.}
    \label{fig:sub-b}
  \end{subfigure}\hfill
  \begin{subfigure}[t]{0.32\textwidth}
    \centering
    \includegraphics[width=\linewidth]{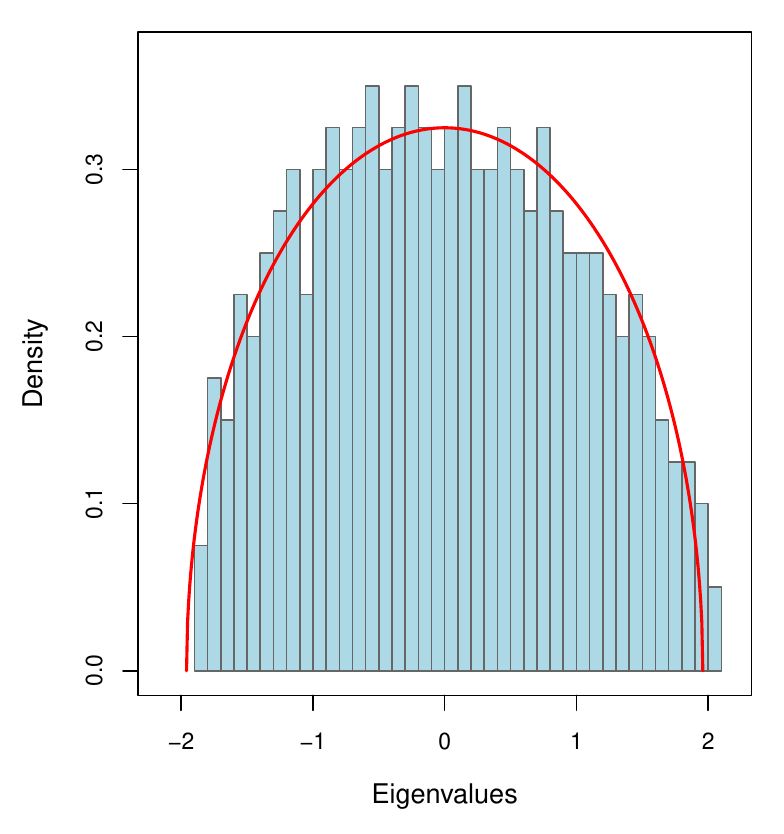}
    \caption{ESD of $\sqrt{n}\Big(\widehat{\fR}_n^{\tau^*}-\fI_p\Big)$.}
    \label{fig:sub-c}
  \end{subfigure}

  \caption{Sample eigenvalue histograms of $\sqrt{n}\Big(\widehat{\fR}_n^D-\fI_p\Big)$, $\sqrt{n}\Big(\widehat{\fR}_n^R-\fI_p\Big)$ and $\sqrt{n}\Big(\widehat{\fR}_n^{\tau^*}-\fI_p\Big)$ with $p=400$, $n=300$. The solid curves are the corresponding limiting semicircle distributions given in Corollary \ref{cor:crc} with $\gamma=p/n=4/3$.}
  \label{fig:three-ESD-figures}
\end{figure}

Figure~\ref{fig:three-ESD-figures} provides a numerical illustration of Corollary~\ref{cor:crc}.

\section{Main results}\label{sec:main}

This section presents our main general results, of which Corollary~\ref{cor:crc} is a special case. Following \cite{MR4185806}, we formulate these results within a general U-statistic framework. Specifically, for each $j\ne k\in[p]$, define
\begin{align}\label{eq-def-Ujk}
\widehat{U}_{jk}
=
\binom{n}{m}^{-1}
\sum_{1\le i_1<i_2<\cdots<i_m\le n}
h\Big(
\binom{X_{i_1,j}}{X_{i_1,k}},
\ldots,
\binom{X_{i_m,j}}{X_{i_m,k}}
\Big),
\end{align}
so that the correlation matrix in \eqref{eq:cons-rank-corr} may equivalently be written as
\begin{align}\label{eq-def-Rjk}
\widehat{R}_{jk}
=
\begin{cases}
1, & \text{if } j=k,\\[0.3em]
\widehat{U}_{jk}, & \text{if } j\ne k.
\end{cases}
\end{align}

To state our general results, we first introduce several basic notions concerning U-statistics.

\begin{definition}[Symmetric kernel]
A kernel $h:(\R^2)^m\to\R$ is said to be {\it symmetric} if
\[
h(\mz_1,\ldots,\mz_m)=h(\mz_{\sigma(1)},\ldots,\mz_{\sigma(m)})
\]
for every permutation $\sigma\in\ml{P}_m$ and every $\mz_1,\ldots,\mz_m\in\R^2$.
\end{definition}

\begin{definition}[Rank-based kernel]\label{def:rank-based}
A kernel $h:(\R^2)^m\to\R$ is said to be {\it rank-based} if, for any vectors $\mz_1,\ldots,\mz_m\in\R^2$, we have
\[
h(\mz_1,\ldots,\mz_m)=h(\mr_1,\ldots,\mr_m),
\]
where $\mr_i=(r_{i,1},r_{i,2})$ records the ranks of $z_{i,1}$ and $z_{i,2}$ among $z_{1,1},\ldots,z_{m,1}$ and $z_{1,2},\ldots,z_{m,2}$, respectively.
\end{definition}

\begin{definition}[$\P_{\mZ}$-degenerate kernel]
Let $\mZ,\mZ_1,\ldots,\mZ_m\in\R^2$ be independent random vectors with common distribution $\P_{\mZ}$. A symmetric kernel $h:(\R^2)^m\to\R$ is said to be {\it $\P_{\mZ}$-degenerate}, or simply {\it degenerate}, if
\[
\var\big(h_1(\mZ;\P_{\mZ})\big)=0.
\]
Here, for each $\ell\in[m]$, we define
\[
h_\ell(\mz_1,\ldots,\mz_\ell;\P_{\mZ})
:=
\E\Big[
h(\mz_1,\ldots,\mz_\ell,\mZ_{\ell+1},\ldots,\mZ_m)
\Big].
\]
\end{definition}

With these notions in place, we are ready to state the main assumptions on the U-statistics under consideration. These assumptions are collected in the following condition. Throughout, let $\P_0$ denote the uniform distribution on $[0,1]$, and let ``$\otimes$'' denote the product of measures. Then $\P_0 \otimes \P_0$ denotes the uniform distribution on $[0,1]^2$.

\begin{assumption}\label{ass:kernel}
Let $\mZ_1,\ldots,\mZ_m\in\R^2$ be independent random vectors sampled from $\P_0\otimes\P_0$. The kernel $h(\cdot)$ is assumed to satisfy the following conditions.
\begin{enumerate}[label=(\roman*)]
\item The kernel $h$ is symmetric, rank-based, and bounded. Moreover, it is mean-zero and degenerate under independent continuous margins, in the sense that
\[
\E\big[h(\mZ_1,\ldots,\mZ_m)\big]=0
\qquad\text{and}\qquad
\E\Big[h_1\big(\mZ_1;\P_0\otimes\P_0\big)^2\Big]=0.
\]

\item There exists a function $g:\R^2\to\R$ such that, with $\mZ_1=(Z_{1,1},Z_{1,2})^\tp$ and $\mZ_2=(Z_{2,1},Z_{2,2})^\tp$,
\[
h_2(\mZ_1,\mZ_2;\P_0\otimes\P_0)
=
g(Z_{1,1},Z_{2,1})\,g(Z_{1,2},Z_{2,2})
\qquad \text{almost surely}.
\]

\item Furthermore, $g(Z_{1,1},Z_{2,1})$ admits the expansion
\[
g(Z_{1,1},Z_{2,1})
=
\sum_{r=1}^\infty \lambda_r \psi_r(Z_{1,1})\psi_r(Z_{2,1}),
\]
where $\{\lambda_r\}$ and $\{\psi_r\}$ are the eigenvalues and eigenfunctions associated with the integral equation
\[
\E\big[g(z_{1,1},Z_{1,2})\psi(Z_{1,2})\big]
=
\lambda\,\psi(z_{1,1}),
\qquad \text{for all } z_{1,1}\in\R,
\]
and satisfy
\[
\lambda_1\ge \lambda_2\ge \cdots \ge 0,
\qquad
\sum_{r=1}^\infty |\lambda_r|\in(0,\infty),
\qquad
\sup_r \|\psi_r\|_\infty<\infty.
\]
\end{enumerate}
\end{assumption}

\begin{remark}
Assumption~\ref{ass:kernel}(i) is exactly Assumption~2.1(i)--(ii) in \cite{MR4185806}. Assumption~\ref{ass:kernel}(ii)--(iii) is stronger than Assumption~2.1(iii) in \cite{MR4185806}. Indeed, under Assumption~\ref{ass:kernel}(ii)--(iii), the kernel $h_2(\mZ_1,\mZ_2;\P_0\otimes\P_0)$ admits the expansion
\[
h_2(\mZ_1,\mZ_2;\P_0\otimes\P_0)
=
\sum_{r,s=1}^\infty \wt{\lambda}_{r,s}\phi_{r,s}(\mZ_1)\phi_{r,s}(\mZ_2),
\]
with eigenvalues $\{\wt{\lambda}_{r,s}\}$ and eigenfunctions $\{\phi_{r,s}\}$ given by
\[
\wt{\lambda}_{r,s}=\lambda_r\lambda_s,
\qquad
\phi_{r,s}(\mZ_1)=\psi_r(Z_{1,1})\psi_s(Z_{1,2}),
\]
and satisfying
\[
\sum_{r,s=1}^\infty \wt{\lambda}_{r,s}
=
\Big(\sum_{r=1}^\infty \lambda_r\Big)^2
<\infty.
\]
Therefore, Assumption~\ref{ass:kernel}(ii)--(iii) is a special case of \citet[Assumption~2.1(iii)]{MR4185806}.
\end{remark}

Based on the work of \cite{MR4185806} and our own calculation, we can prove the following proposition. It shows that the three consistent rank correlations mentioned in Section~\ref{sec:crc} all satisfy Assumption~\ref{ass:kernel}. 

\begin{proposition}\label{prop-kernel-properties}
The following statements hold.
\begin{enumerate}[label=(\roman*)]
\item {\bf Hoeffding's $D$.} The kernel $h_D$ in \eqref{eq:hoeffding} satisfies Assumption~\ref{ass:kernel}, with the corresponding function $g(\cdot)$ and associated eigenvalues and eigenfunctions given by
\[
g(x,y)=\frac{\sqrt{3}}{6}\Big(3x^2+3y^2-6\max(x,y)+2\Big),
\qquad
\lambda_r=\frac{\sqrt{3}}{\pi^2 r^2},
\qquad
\psi_r(x)=\sqrt{2}\cos(\pi r x).
\]

\item {\bf Blum--Kiefer--Rosenblatt's $R$.} The kernel $h_R$ in \eqref{eq:bkr} satisfies Assumption~\ref{ass:kernel}, with
\[
h_{R,2}=2h_{D,2}.
\]

\item {\bf Bergsma--Dassios--Yanagimoto's $\tau^*$.} The kernel $h_{\tau^*}$ in \eqref{eq:bdy} satisfies Assumption~\ref{ass:kernel}, with
\[
h_{\tau^*,2}=3h_{D,2}.
\]
\end{enumerate}
\end{proposition}

Moreover, the three consistent rank correlations discussed in Section~\ref{sec:crc} are all generated by kernels that are consistent for testing independence.

\begin{proposition}[Proposition~3 in \cite{shi2020power}]
Let $\mZ=(Z_1,Z_2)^\tp\in\R^2$ be Lebesgue absolutely continuous, and let $\mZ_1,\mZ_2,\ldots$ be independent copies of $\mZ$. Then the following four statements are equivalent:

\begin{enumerate}
  \item $\E\big[h_D(\mZ_1,\ldots,\mZ_5)\big]=0$;
  \item $\E\big[h_R(\mZ_1,\ldots,\mZ_6)\big]=0$;
  \item $\E\big[h_{\tau^*}(\mZ_1,\ldots,\mZ_4)\big]=0$;
  \item $Z_1$ is independent of  $Z_2$.
\end{enumerate}

\end{proposition}

We are now in a position to state the main theorem of this paper. To this end, define the correlation matrix $\widehat{\fR}_n=(\widehat{R}_{jk})$ with the entries $\widehat{R}_{jk}$ given in \eqref{eq-def-Rjk}, and its standardized version
\[
\widehat{\fW}_n=\sqrt{n}\,(\widehat{\fR}_n-\fI_p).
\]

\begin{thm}\label{thm:sc-law-main}
Assume \eqref{eq:assumption}, \eqref{eq:asymptotics}, and Assumption~\ref{ass:kernel}. Then
\[
F^{\widehat{\fW}_n}\Rightarrow \operatorname{W}(\vartheta)
\qquad \text{almost surely},
\]
with the radius parameter
\[
\vartheta:=m(m-1)\sqrt{2\gamma}\,\sum_{r=1}^\infty \lambda_r^2
\]
which is a well-defined positive constant.
\end{thm}


\section{Proof of main results}

We first introduce some additional notation. Write
\[
\R^+:=\{u\in\R:u>0\}, \qquad
\bC^+:=\{u+\mathrm{i}v:u,v\in\bR,\ v>0\}, \qquad
\text{and} \qquad
\bN=\{1,2,3,\ldots\},
\]
 where $\mathrm{i}=\sqrt{-1}$ is the imaginary unit. For any complex number $z$, let $\Im z$ denote its imaginary part. For two random variables $X$ and $Y$, we write $X\perp Y$ if $X$ and $Y$ are independent. For two sequences $\{a_n\}_{n=1}^\infty,\{b_n\}_{n=1}^\infty\subset\R^+$, we write $a_n\lesssim b_n$ if $a_n\leq Cb_n$ for some constant $C$ and all sufficiently large $n$, $a_n\sim b_n$ if $a_n=b_n(1+o(1))$, and $a_n\asymp b_n$ if both $a_n\lesssim b_n$ and $b_n\lesssim a_n$ hold. 

For a matrix $\fX\in\bC^{m\times n}$, let $\tr(\fX)$ denote its trace and $\|\fX\|_{S_p}$ denote its Schatten $p$-norm defined by
$$
\|\fX\|_{S_p}:=\Bigl(\sum_j \sigma_j(\fX)^p\Bigr)^{1/p},
\qquad p\geq 1,
$$
where $\{\sigma_j(\fX)\}_j$ are the singular values of $\fX$. In particular, 
\[
\|\fX\|_{\mathrm{F}}:=\|\fX\|_{S_2}~~~{\rm and}~~~ \|\fX\|_{\mathrm{op}}:=\|\fX\|_{S_\infty} 
\]
correspond to the Frobenius and spectral norms, respectively. For a square matrix $\fA$, write
\[
D_0(\fA):=\fA-\Diag(\fA).
\]
Here $\Diag(\fA)$ denotes the diagonal matrix made with the diagonal entries of $\fA$.
For a probability measure $\mu$ on $\bR$, its Stieltjes transform is defined by
\[
s(z):=\int_\bR \frac{\mu(dx)}{x-z},
\qquad z\in\bC^+.
\]
In particular, the Stieltjes transform of Wigner's semicircle distribution $\mathrm{W}(r)$ is
\be\label{eq-Stieljes-transform-W(r)}
m_r(z)=-\frac{2}{r^2}\Bigl(z-\sqrt{z^2-r^2}\Bigr),
\qquad z\in\bC^+.
\ee
For a symmetric matrix $\fA\in\bR^{p\times p}$, we write $s_\fA(z)$ for the Stieltjes transform of $F^{\fA}$.

\subsection{A proof outline}\label{section:proof-outline}

In our proof, without loss of generality, we assume that the entries $\{X_{ij}: i\in[n],\, j\in[p]\}$ are independently distributed according to $\P_0$. This does not affect the distribution of $\widehat{\fR}_n$ under Assumption~\ref{eq:assumption}, since each $\widehat{R}_{jk}$ is rank-based.

We first reduce the analysis of $\widehat{\fR}_n$ to that of its leading term. Since the U-statistic $\widehat{R}_{jk}$ is degenerate under Assumption~\ref{ass:kernel}, the first nonvanishing term in its Hoeffding decomposition is a second-order U-statistic, given by
\begin{align}
  \widehat{R}_{jk}^{(1)}
  :=\frac{m(m-1)}{n(n-1)}\sum_{1\leq i_1<i_2\leq n}
  h_{2}\l(\binom{X_{i_1,j}}{X_{i_1,k}},\binom{X_{i_2,j}}{X_{i_2,k}}\r).
  \nonumber
\end{align}
By Assumption~\ref{ass:kernel},
\begin{align}
  \widehat{R}_{jk}^{(1)}
  &=\frac{m(m-1)}{n(n-1)}\sum_{1\leq i_1<i_2\leq n}
  g(X_{i_1,j},X_{i_2,j})g(X_{i_1,k},X_{i_2,k})
  \nonumber \\
  &=\frac{m(m-1)}{n(n-1)}\sum_{1\leq i_1<i_2\leq n}\sum_{r,s=1}^\infty
  \lambda_r\lambda_s\psi_r(X_{i_1,j})\psi_s(X_{i_1,k})\psi_r(X_{i_2,j})\psi_s(X_{i_2,k}).
  \label{eq-first-term-Hoeffding}
\end{align}
Setting $\widehat{R}_{jj}^{(1)}:=1$, we define the leading term of $\widehat{\fR}_n$ by
$$
\widehat{\fR}_n^{(1)}:=[\widehat{R}_{jk}^{(1)}]_{j,k\in [p]}.
$$
Later we will show that, at the LSD scale, the difference between $\widehat{\fR}_n^{(1)}$ and $\widehat{\fR}_n$ is negligible. We may therefore restrict attention to the LSD of $\widehat{\fR}_n^{(1)}$.

Next, following the analysis in \cite{shi2020power}, we apply a truncation argument to handle the infinite sum over $r$ and $s$ in \eqref{eq-first-term-Hoeffding}. More precisely, for $T\in\N$, define the truncated coefficient
$$
\widehat{R}_{jk;T}^{(1)}
:=
\frac{m(m-1)}{n(n-1)}
\sum_{1\leq i_1<i_2\leq n}\sum_{r,s=1}^T
\lambda_r\lambda_s
\psi_r(X_{i_1,j})\psi_s(X_{i_1,k})
\psi_r(X_{i_2,j})\psi_s(X_{i_2,k})
$$
for $j\neq k\in[p]$, and set
\[
\widehat{R}_{jj;T}^{(1)}:=1
\]
for $j\in[p]$. We then introduce the truncated leading matrix
\[
\widehat{\fR}_{n;T}^{(1)}:=\bigl(\widehat{R}_{jk;T}^{(1)}\bigr)_{j,k\in[p]}.
\]
Later we will show that, at the LSD scale, for all sufficiently large $T$, the difference between $\widehat{\fR}_{n;T}^{(1)}$ and $\widehat{\fR}_{n}^{(1)}$ is negligible. It therefore remains to analyze the LSD of $\widehat{\fR}_{n;T}^{(1)}$.


Next, we rewrite $\widehat{R}_{jk;T}^{(1)}$ as
$$
\widehat{R}_{jk;T}^{(1)}
=
\frac{m(m-1)}{n(n-1)}
\sum_{1\leq i_1<i_2\leq n}
\Bigl(\sum_{r=1}^T\lambda_r\psi_r(X_{i_1,j})\psi_r(X_{i_2,j})\Bigr)
\Bigl(\sum_{s=1}^T\lambda_s\psi_s(X_{i_1,k})\psi_s(X_{i_2,k})\Bigr).
$$
This representation shows that the off-diagonal part of $\widehat{\fR}_{n;T}^{(1)}$ admits a high-dimensional Gram matrix representation. To formalize this, define
\[
\cM_2:=\{(i_1,i_2): 1\leq i_1<i_2\leq n\},
\qquad
M:=|\cM_2|=\frac{n(n-1)}{2},
\]
and set
$$
A_{\mi,k}:=\sum_{r=1}^T \lambda_r\psi_r(X_{i_1,k})\psi_r(X_{i_2,k}),
\qquad
\text{for all }\mi=(i_1,i_2)\in \cM_2 \text{ and } k\in[p].
$$
Then, for $j\neq k\in[p]$, we have
$$
\widehat{R}_{jk;T}^{(1)}
=
\frac{m(m-1)}{n(n-1)}\,\mA_{\cdot,k}^\tp \mA_{\cdot,j},
$$
where $\mA_{\cdot,k}:=(A_{\mi,k})_{\mi\in\cM_2}\in\R^M$. Introducing
\[
\fA:=(A_{\mi,k})_{\mi\in\cM_2,\,k\in[p]},
\]
and recalling that $\widehat{R}_{jj;T}^{(1)}:=1$, the preceding identity yields the off-diagonal Gram representation
$$
D_0\bigl(\widehat{\fR}_{n;T}^{(1)}\bigr)
=
\frac{m(m-1)}{n(n-1)}\,D_0\bigl(\fA^\tp\fA\bigr).
$$
Equivalently, $\widehat{\fR}_{n;T}^{(1)}$ is obtained from
$\frac{m(m-1)}{n(n-1)}\fA^\tp\fA$
by replacing its diagonal entries with ones. Here the matrix $\fA$ depends on both $n$ and $T$. Thus, the off-diagonal part of $\widehat{\fR}_{n;T}^{(1)}$ admits an explicit Gram matrix representation. Accordingly, the associated rescaled matrix can be written as the following diagonal-removed Gram matrix:
\begin{align*}
\widehat{\fW}_{n;T}^{(1)}
&:=\sqrt{n}\bigl(\widehat{\fR}_{n;T}^{(1)}-\fI_p\bigr)\\
&=\frac{m(m-1)}{\sqrt{n}(n-1)}\,D_0\bigl(\fA^\tp\fA\bigr)\\
&=\frac{m(m-1)}{\sqrt{n}(n-1)}
\bigl(\fA^\tp\fA-\Diag(\fA^\tp\fA)\bigr).
\end{align*}

As will be shown later, the LSD of $\widehat{\fW}_{n}$ is well approximated by that of $\widehat{\fW}_{n;T}^{(1)}$. Moreover, the matrix $\fA$ appearing in $\widehat{\fW}_{n;T}^{(1)}$ satisfies the following properties:
\begin{enumerate}[label=(\roman*)]
\item The columns of $\fA$ are independent and identically distributed (i.i.d.), and each column has covariance matrix $\sigma_T^2\fI_M$ for some constant $\sigma_T>0$ that possibly depends on $T$.
\item Under \eqref{eq:asymptotics}, the dimension of $\fA$ is of order $O(n^2)\times n$.
\end{enumerate}
These observations strongly suggest that $\widehat{\fW}_{n;T}^{(1)}$, and hence also $\widehat{\fW}_n$, should satisfy the semicircle law, in view of the seminal result of \cite{BaiYin1988SC} on the LSD of Gram matrices in the ultra-high-dimensional regime.

Unfortunately, several technical obstacles remain. In particular, the results of \cite{BaiYin1988SC} are not directly applicable here, since the entries of $\fA$ are correlated. Establishing the semicircle law therefore requires additional control of the concentration of certain random quadratic forms associated with the columns of $\fA$. This issue will be addressed in Section \ref{sec::proof-of-prop-EsGz-convergence}, Lemma \ref{Lemma-Quadratic-form-control}.

To summarize, our proof proceeds through two approximation steps for the original sample correlation matrix $\widehat{\fR}_n$. First, we approximate $\widehat{\fR}_n$ entrywise by a second-order U-statistic, thereby obtaining $\widehat{\fR}_n^{(1)}$. Second, we truncate the infinite sum in \eqref{eq-first-term-Hoeffding} to obtain $\widehat{\fR}_{n;T}^{(1)}$. After these two approximations, the off-diagonal part of $\widehat{\fR}_{n;T}^{(1)}$ admits the useful Gram matrix representation
$$
D_0\bigl(\widehat{\fR}_{n;T}^{(1)}\bigr)
=
\frac{m(m-1)}{n(n-1)}D_0\bigl(\fA^\tp\fA\bigr),
$$
although the entries of $\fA$ remain correlated.

In Section~\ref{section-sc-limit-Wnt1}, we establish the semicircle limit of
$\widehat{\fW}_{n;T}^{(1)}=\sqrt{n}\bigl(\widehat{\fR}_{n;T}^{(1)}-\fI_p\bigr)$,
which serves as a good approximation to the target matrix $\widehat{\fW}_n$, by following the classical Stieltjes transform method of \cite{bai2010spectral}. The main technical challenge is to control the dependence among the entries of $\fA$. We then prove the LSD of $\widehat{\fW}_n$ in Section~\ref{section-main-thm-pf} by carefully quantifying the errors introduced in these two approximation steps.

\subsection{Technical preparation: Semicircle law of $\widehat{\fW}_{n;T}^{(1)}$}\label{section-sc-limit-Wnt1}

We retain the notation and definitions introduced in Section~\ref{section:proof-outline}. In addition, for $\mi,\mj\in \cM_2$, with
\[
\mi=(i_1,i_2)
\qquad \text{and} \qquad
\mj=(j_1,j_2),
\]
we
define
\[
\delta_{\mu\nu}:=\ind(\mu=\nu).
\]
We further write
\[
\mi\cap \mj\neq \emptyset
\qquad \text{if} \qquad
\{i_1,i_2\}\cap \{j_1,j_2\}\neq \emptyset.
\]

Since the dependence among the entries of $\fA$ is the central issue in the proof, we begin by identifying the correlation structure of $\fA$. Assumption~\ref{ass:kernel} implies that for all $i\in[n]$ and $k\in[p]$,
\be\label{eq-orthogonal-psi}
\bE\big[\psi_r(X_{ik})\big]=0,
\qquad
\bE\big[\psi_r(X_{ik})\psi_s(X_{ik})\big]=\delta_{rs}.
\ee
Consequently, the entrywise variance of $\fA$ is given by
\[
\E\big[A_{\mi,k}^2\big]=\sum_{r=1}^T \lambda_r^2.
\]
Accordingly, we normalize $\fA$ by introducing
\[
\sigma_T:=\Bigl(\sum_{r=1}^T\lambda_r^2\Bigr)^{1/2}
\qquad \text{and} \qquad
\wt{\fA}:=\sigma_T^{-1}\fA,
\]
so that each entry of $\wt{\fA}$ has unit variance. Note that
$\wt{\fA}=(\wt{A}_{\mi,k})_{\mi\in\cM_2,\ k\in[p]}$, where
$$
\wt{A}_{\mi,k}
=
\sigma_T^{-1}\sum_{r=1}^T \lambda_r\psi_r(X_{i_1,k})\psi_r(X_{i_2,k}),
\qquad
\mi=(i_1,i_2)\in\cM_2.
$$
The following five properties of $\wt{\fA}$ are then immediate.

\begin{enumerate}[label=(\roman*)]
\item The column vectors $\{\wt{\mA}_{\cdot,k}: k\in[p]\}$ are i.i.d.\ random vectors in $\bR^M$.
\item By Assumption~\ref{ass:kernel}, the entries of $\wt{\fA}$ are uniformly bounded:
$$
\sup_{\mi\in\cM_2,\ k\in[p]} |\wt{A}_{\mi,k}|
\leq
\sigma_T^{-1}\Bigl(\sup_r \|\psi_r\|_\infty\Bigr)^2 \sum_{r=1}^\infty |\lambda_r|
<\infty.
$$
\item The random vector $\wt{\mA}_{\cdot,1}$ is isotropic in $\R^M$, since
$$
\E[\wt{A}_{\mi,1}]=0,
\qquad
\E[\wt{A}_{\mi,1}^2]=1,
\qquad
\E[\wt{A}_{\mi,1}\wt{A}_{\mj,1}]=\delta_{\mi\mj}.
$$
\item If $\mi\cap\mj=\emptyset$, then the random variables $\wt{A}_{\mi,1}$ and $\wt{A}_{\mj,1}$ are independent.
\item The collection $\{\wt{A}_{\mi,1}:\mi\in\cM_2\}$ is not jointly independent. For example, let $\mi_1=(1,2)$, $\mi_2=(2,3)$, and $\mi_3=(1,3)$. Using \eqref{eq-orthogonal-psi} once again,
$$
\mathbb{E}[\wt{A}_{\mi_1,1}\wt{A}_{\mi_2,1}\wt{A}_{\mi_3,1}]
=
\sigma_T^{-3}\sum_{r=1}^T \lambda_r^3
\neq
\bE[\wt{A}_{\mi_1,1}]\,\bE[\wt{A}_{\mi_2,1}]\,\bE[\wt{A}_{\mi_3,1}].
$$
\end{enumerate}

Nevertheless, the dependence within $(\wt{A}_{\mi,k})_{\mi\in\cM_2}$ is relatively weak. The following lemma captures this weak dependence by characterizing the vanishing of certain cross moments, and it will play a key role in the proof.

\begin{lem}\label{lemma-cross-moments}
Consider arbitrary $\mi_1,\ldots,\mi_Q\in\cM_2$, where $\mi_q=(i_1^q,i_2^q)$ for each $q\in[Q]$. If, among the $2Q$ indices
\[
i_1^1,i_2^1,\ldots,i_1^Q,i_2^Q,
\]
there exists an index that appears exactly once, then
\[
\E\Bigl[\prod_{q=1}^Q \wt{A}_{\mi_q,1}\Bigr]=0.
\]
\end{lem}

We are now in a position to formulate the semicircle law for $\widehat{\fW}_{n;T}^{(1)}$. Define
$$
\fG_{n;T}
:=
\sqrt{\frac{M}{p}}
\Bigl(
\frac{1}{M}\wt{\fA}^\tp\wt{\fA}
-
\Diag\Bigl(\frac{1}{M}\wt{\fA}^\tp\wt{\fA}\Bigr)
\Bigr).
$$
Then
\bea\label{eq-hat-W-n-T-(1)}
\widehat{\fW}_{n;T}^{(1)}
=
\frac{m(m-1)\sqrt{pM}\sigma_T^2}{\sqrt{n}(n-1)}\,\fG_{n;T},
\eea
where the matrix $\fG_{n;T}$ has a structure analogous to that of the ultra-high-dimensional Gram matrices appearing in \cite{BaiYin1988SC}, which exhibit semicircle limits. Lemma~\ref{lemma-cross-moments} shows that the dependence within $(\wt{A}_{\mi,k})_{\mi\in\cM_2}$ is sufficiently weak for $\fG_{n;T}$ to retain a semicircle limit. This, in turn, yields the semicircle limit of $\widehat{\fW}_{n;T}^{(1)}$ up to a deterministic scaling factor.

The semicircle law for $\widehat{\fW}_{n;T}^{(1)}$ is stated in the following propositions.

\begin{proposition}\label{prop-stieljes-concentration}
There exists a constant $C>0$ such that for any $\varepsilon>0$, $z=u+\mathrm{i}v\in\bC^+$, and $T\in\N$,
$$
\bP\Bigl(\bigl|s_{\widehat{\fW}_{n;T}^{(1)}}(z)-\E\big[s_{\widehat{\fW}_{n;T}^{(1)}}(z)\big]\bigr|>\varepsilon\Bigr)
\leq
\frac{C}{\varepsilon^4 v^4 p^2}.
$$
\end{proposition}

\begin{proposition}\label{prop-EsG(z)-convergence}
Under \eqref{eq:asymptotics}, for any fixed $T\in\N$ and $z\in\bC^+$,
\[
\lim_{n\to\infty}\E\Big[s_{\widehat{\fW}_{n;T}^{(1)}}(z)\Big]=m_{r_T}(z),
\]
where $\vartheta=\lim_{T\to\infty} r_T$  with 
\[
r_T:=m(m-1)\sqrt{2\gamma}\sum_{r=1}^T \lambda_r^2,
\]
and $m_{r_T}(z)$ denotes the Stieltjes transform of $\operatorname{W}(r_T)$ given in \eqref{eq-Stieljes-transform-W(r)}.
\end{proposition}


\subsection{Proof of Theorem \ref{thm:sc-law-main}}\label{section-main-thm-pf}

Fix $\varepsilon>0$, and for $z\in\bC^+$, write $v=\Im z>0$. Let $r_T$ be as defined in Proposition~\ref{prop-EsG(z)-convergence}. Since $\lim_{T\to\infty} r_T=\vartheta$, there exists a constant $T_\varepsilon\in\mathbb{N}$ such that
\bea\label{eq-5steps-bound-0}
\big|m_{r_T}(z)-m_{\vartheta}(z)\big|<\frac{\varepsilon}{5},
\qquad \text{for all } T>T_\varepsilon.
\eea
Applying the same argument as in the proof of Proposition~\ref{prop-stieljes-concentration} to
\[
\widehat{\fW}_n^{(1)}:=\sqrt{n}\bigl(\widehat{\fR}_n^{(1)}-\fI_p\bigr),
\]
we obtain that there exists a constant $C>0$ such that
\bea\label{eq-5steps-bound-1}
\bP\Bigl(\big|s_{\widehat{\fW}_{n}^{(1)}}(z)-\E\big[s_{\widehat{\fW}_{n}^{(1)}}(z)\big]\big|>\frac{\varepsilon}{5}\Bigr)
\leq
\frac{C}{\varepsilon^4 v^4 p^2}.
\eea
By the variance formula for U-statistics \citep[p.~189]{serfling}, as $n\to\infty$,
$$
\E\Bigl[\big|\widehat{R}_{12}^{(1)}-\widehat{R}_{12;T}^{(1)}\big|^2\Bigr]
=
\frac{\xi_T}{n^2}
+
O\Bigl(\frac{1}{n^3}\Bigr)
$$
for some positive constant $\xi_T$ depending on $T$ such that $\xi_T\to 0$ as $T\to\infty$.
We then choose $\wt{T}_{\varepsilon}\in\mathbb{N}$ such that
$$
\wt{T}_{\varepsilon}>T_{\varepsilon}
\qquad \text{and} \qquad
\xi_{\wt{T}_{\varepsilon}}<\frac{\varepsilon^2 v^4}{100\gamma}.
$$

By exchangeability of the samples and coordinates, all off-diagonal entries of
$\widehat{\fR}_{n}^{(1)}-\widehat{\fR}_{n;\wt{T}_{\varepsilon}}^{(1)}$
have the same distribution, and likewise all off-diagonal entries of
$\widehat{\fR}_n-\widehat{\fR}_n^{(1)}$
have the same distribution. We shall use this symmetry in the Frobenius norm bounds below.
    
Then, by Lemma~\ref{lem-perturbation-frob-norm} below, for all sufficiently large $n$,
\bea\label{eq-5steps-bound-2}
\big|\E\big[s_{\widehat{\fW}_{n}^{(1)}}(z)\big]-\E\big[s_{\widehat{\fW}_{n;\wt{T}_{\varepsilon}}^{(1)}}(z)\big]\big|
&\leq
\sqrt{\E\Big[\big|s_{\widehat{\fW}_{n}^{(1)}}(z)-s_{\widehat{\fW}_{n;\wt{T}_{\varepsilon}}^{(1)}}(z)\big|^2\Big]}\\
&\leq
\sqrt{\frac{n}{pv^4}\E\Big[
\big\|\widehat{\fR}_{n}^{(1)}-\widehat{\fR}_{n;\wt{T}_{\varepsilon}}^{(1)}\big\|_{\operatorname{F}}^2
\Big]}\\
&\leq
\sqrt{\frac{n}{pv^4}\cdot p(p-1)\cdot \frac{2\xi_{\wt{T}_{\varepsilon}}}{n^2}}\\
&\leq
\frac{\varepsilon}{5}.
\eea
By Proposition~\ref{prop-EsG(z)-convergence} and \eqref{eq-hat-W-n-T-(1)},
$$
\lim_{n\to\infty}\E\big[s_{\widehat{\fW}_{n;\wt{T}_{\varepsilon}}^{(1)}}(z)\big]
=
m_{r_{\wt{T}_{\varepsilon}}}(z).
$$
Hence, for all sufficiently large $n$,
\bea\label{eq-5steps-bound-3}
\big|\E\big[s_{\widehat{\fW}_{n;\wt{T}_{\varepsilon}}^{(1)}}(z)\big]-m_{r_{\wt{T}_{\varepsilon}}}(z)\big|
<
\frac{\varepsilon}{5}.
\eea
In addition, by \citet[Equation~4, p.~190]{serfling}, there exists a constant $D>0$ such that
$$
\E\Big[\big|\widehat{R}_{12}-\widehat{R}_{12}^{(1)}\big|^4\Big]\leq \frac{D}{n^6}.
$$
Then, by Lemma~\ref{lem-perturbation-frob-norm} and the Cauchy--Schwarz inequality,
\bea\label{eq-5steps-bound-4}
\P\Big(\big|s_{\widehat{\fW}_n}(z)-s_{\widehat{\fW}_n^{(1)}}(z)\big|>\frac{\varepsilon}{5}\Big)
&\leq
\frac{C'}{\varepsilon^4 p^2 v^8}
\E\Big[\big\|\widehat{\fW}_n-\widehat{\fW}_n^{(1)}\big\|_{\operatorname{F}}^4\Big]\\
&\leq
\frac{C'n^2}{\varepsilon^4 p^2 v^8}
\sum_{j\neq k=1}^p\sum_{j'\neq k'=1}^p
\E\Big[
\big|\widehat{R}_{jk}-\widehat{R}_{jk}^{(1)}\big|^2
\big|\widehat{R}_{j'k'}-\widehat{R}_{j'k'}^{(1)}\big|^2
\Big]\\
&\leq
\frac{C'p^2}{\varepsilon^4 n^4 v^8}.
\eea
Combining \eqref{eq-5steps-bound-0}--\eqref{eq-5steps-bound-4}, we obtain
$$
\P\Big(\big|s_{\widehat{\fW}_{n}}(z)-m_{\vartheta}(z)\big|>\varepsilon\Big)
\leq
\frac{C}{\varepsilon^4 v^4 p^2}
+
\frac{C'p^2}{\varepsilon^4 n^4 v^8}.
$$
By the Borel--Cantelli lemma and the assumption $p\sim n$, it follows that for any fixed $z\in\bC^+$,
$$
s_{\widehat{\fW}_{n}}(z)-m_{\vartheta}(z)\to 0,
\qquad \text{almost surely.}
$$
The proof is therefore completed by a standard argument based on Vitali's convergence theorem \citep[Lemma~2.14]{bai2010spectral} and the continuity theorem for Stieltjes transforms \citep[Theorem~B.9]{bai2010spectral}.

\section{Proof of auxiliary results}

\subsection{Proof of Propositions}
This section proves the propositions stated in the main text. In this process, we also develop some auxiliary lemmas, whose proofs are put in Section \ref{section-proof-of-lemmas}. 

\subsubsection{Proof of Proposition \ref{prop-kernel-properties}}

We first focus on Part~(i). To this end, by Example~2.1 of \cite{MR4185806}, we have
$$
h_2\!\left(
\begin{pmatrix}
z_{1,1}\\
z_{1,2}
\end{pmatrix},
\begin{pmatrix}
z_{2,1}\\
z_{2,2}
\end{pmatrix}
\right)
=
\sum_{i,j=1}^\infty \frac{6}{\pi^4 i^2 j^2}
\cos(\pi i z_{1,1})\cos(\pi i z_{2,1})
\cos(\pi j z_{1,2})\cos(\pi j z_{2,2}).
$$
Hence $h_2$ can be written as
$$
h_2\!\left(
\begin{pmatrix}
z_{1,1}\\
z_{1,2}
\end{pmatrix},
\begin{pmatrix}
z_{2,1}\\
z_{2,2}
\end{pmatrix}
\right)
=
\sum_{r,s=1}^{\infty}
\lambda_r\lambda_s
\psi_r(z_{1,1})\psi_r(z_{2,1})
\psi_s(z_{1,2})\psi_s(z_{2,2}),
$$
where
\[
\lambda_r=\frac{\sqrt{3}}{\pi^2 r^2},
\qquad
\psi_r(x)=\sqrt{2}\cos(\pi r x).
\]
Equivalently,
$$
h_2\!\left(
\begin{pmatrix}
z_{1,1}\\
z_{1,2}
\end{pmatrix},
\begin{pmatrix}
z_{2,1}\\
z_{2,2}
\end{pmatrix}
\right)
=
\left(
\sum_{r=1}^{\infty}
\lambda_r \psi_r(z_{1,1})\psi_r(z_{2,1})
\right)
\left(
\sum_{s=1}^{\infty}
\lambda_s \psi_s(z_{1,2})\psi_s(z_{2,2})
\right).
$$
Define
\be\label{eq-gxy-proof}
g(x,y)=\sum_{r=1}^{\infty}\lambda_r\psi_r(x)\psi_r(y),
\qquad x,y\in[0,1].
\ee
Then
$$
h_2\!\left(
\begin{pmatrix}
z_{1,1}\\
z_{1,2}
\end{pmatrix},
\begin{pmatrix}
z_{2,1}\\
z_{2,2}
\end{pmatrix}
\right)
=
g(z_{1,1},z_{2,1})\,g(z_{1,2},z_{2,2}).
$$

We next show that $g(x,y)$ is given by the expression stated in Proposition~\ref{prop-kernel-properties}. By \eqref{eq-gxy-proof},
\begin{align*}
g(x,y)
&=
\sum_{r=1}^{\infty}
\frac{\sqrt{3}}{\pi^2 r^2}\cdot 2\cos(\pi r x)\cos(\pi r y)\\
&=
\frac{2\sqrt{3}}{\pi^2}\cdot \frac12
\left(
\sum_{r=1}^{\infty}\frac{\cos(\pi r(x+y))}{r^2}
+
\sum_{r=1}^{\infty}\frac{\cos(\pi r(x-y))}{r^2}
\right).
\end{align*}
We use the identity
$$
\sum_{r=1}^{\infty}\frac{\cos(r\theta)}{r^2}
=
\frac{\pi^2}{6}
-\frac{\pi|\theta|}{2}
+\frac{\theta^2}{4},
\qquad \text{for all }\theta\in[-2\pi,2\pi].
$$
Substituting $\theta=\pi(x+y)$ and $\theta=\pi(x-y)$ yields
\begin{align*}
g(x,y)
&=
\frac{\sqrt{3}}{6}
\left(
2+3x^2+3y^2-3x-3y-3|x-y|
\right)\\
&=
\frac{\sqrt{3}}{6}
\left(
3x^2+3y^2-6\max(x,y)+2
\right).
\end{align*}
This proves Part~(i). Parts~(ii) and~(iii) follow from Examples~2.2 and~2.3 of \cite{MR4185806}, respectively.

\subsubsection{Proof of Proposition \ref{prop-stieljes-concentration}}

Recall that
\[
\widehat{\fW}_{n;T}^{(1)}=b_nD_0(\fA^\tp\fA),
\qquad {\rm with}~~
b_n:=\frac{m(m-1)}{\sqrt{n}(n-1)}.
\]
Define $\rA_0:=\{\emptyset\}$ and, for $1\leq k\leq p$,
\[
\rA_k:=\sigma\Bigl(\{\mA_{\cdot,j}:1\leq j\leq k\}\Bigr),
\qquad
\E_{\rA_k}[\cdot]:=\E[\cdot\mid \rA_k].
\]
Let $\fA_{-k}$ denote the $M\times (p-1)$ matrix obtained by removing the $k$th column $\mA_{\cdot,k}$ from $\fA$. Then
$$
s_{\widehat{\fW}_{n;T}^{(1)}}(z)-\E\big[s_{\widehat{\fW}_{n;T}^{(1)}}(z)\big]
=
\sum_{k=1}^p Q_k,
$$
where
\bea\nonumber
Q_k
&=
(\E_{\rA_k}-\E_{\rA_{k-1}})
\Bigl[
\frac{1}{p}\tr\bigl(b_n D_0(\fA^\tp \fA)-z\fI\bigr)^{-1}
\Bigr]\\
&=
(\E_{\rA_k}-\E_{\rA_{k-1}})
\Bigl[
\frac{1}{p}\tr\bigl(b_n D_0(\fA^\tp \fA)-z\fI\bigr)^{-1}
-
\frac{1}{p}\tr\bigl(b_n D_0(\fA_{-k}^\tp \fA_{-k})-z\fI\bigr)^{-1}
\Bigr].
\eea
By Lemma~\ref{lem-trace-resolvant-minor} and denoting $v=\Im z$, we have
$$
\Bigl|
\tr\bigl(b_nD_0(\fA^\tp \fA)-z\fI\bigr)^{-1}
-
\tr\bigl(b_nD_0(\fA_{-k}^\tp \fA_{-k})-z\fI\bigr)^{-1}
\Bigr|
\leq
\frac{1}{v}.
$$
Therefore,
$$
|Q_k|\leq \frac{2}{pv}.
$$
Applying Lemma~\ref{lem-burkholder}, we obtain
$$
\E\Bigl[\bigl|s_{\widehat{\fW}_{n;T}^{(1)}}(z)-\E[s_{\widehat{\fW}_{n;T}^{(1)}}(z)]\bigr|^4\Bigr]
=
\E\Bigl[\Bigl|\sum_{k=1}^p Q_k\Bigr|^4\Bigr]
\leq
C\,\E\Bigl[\Bigl(\sum_{k=1}^p |Q_k|^2\Bigr)^2\Bigr]
\leq
\frac{16C}{v^4p^2},
$$
where $C>0$ is a constant. The stated probability bound then follows from Markov's inequality.

\subsubsection{Proof of Proposition \ref{prop-EsG(z)-convergence}}\label{sec::proof-of-prop-EsGz-convergence}

By \eqref{eq-hat-W-n-T-(1)}, it suffices to show that
$$
\lim_{n\to\infty}\E[s_{\fG_{n;T}}(z)]=m_2(z),
$$
where
$$
m_2(z)=\frac{-z+\sqrt{z^2-4}}{2}.
$$
Indeed, \eqref{eq-hat-W-n-T-(1)} identifies $\widehat{\fW}_{n;T}^{(1)}$ as a deterministic rescaling of $\fG_{n;T}$. Hence, if the limiting Stieltjes transform of $\fG_{n;T}$ is $m_2(z)$, then the limiting Stieltjes transform of $\widehat{\fW}_{n;T}^{(1)}$ is exactly $m_{r_T}(z)$.

Let
\[
y_n:=\frac{p}{M},
\qquad
\fS_n:=\frac{1}{M}\wt{\fA}^\tp\wt{\fA}.
\]
In this subsection, we write $\fG=\fG_{n;T}$ and $\fS=\fS_n$ for brevity. Then
\bea\label{sGsS-transform}
s_{\fG}(z)=\frac{1}{p}\tr\bigl(y_n^{-1/2}D_0(\fS)-z\fI\bigr)^{-1},
\qquad z\in\bC^+.
\eea
Let $v=\Im z>0$, let
\[
\wt{\mA}_{\cdot,k}:=(\wt{A}_{\mu,k})_{\mu\in\cM_2}\in\R^M,
\]
let $\wt{\fA}_{-k}$ be the $M\times (p-1)$ matrix obtained by removing the $k$th column $\wt{\mA}_{\cdot,k}$ from $\wt{\fA}$, and define
\[
\fS_{-k}:=\frac{1}{M}\wt{\fA}_{-k}^\tp\wt{\fA}_{-k}.
\]
By \citet[Theorem~A.4]{bai2010spectral},
$$
s_{\fG}(z)
=
\frac{1}{p}\sum_{k=1}^p
\frac{1}{
-z-y_n^{-1}M^{-2}\wt{\mA}_{\cdot,k}^\tp
\wt{\fA}_{-k}
\bigl(y_n^{-1/2}D_0(\fS_{-k})-z\bigr)^{-1}
\wt{\fA}_{-k}^\tp \wt{\mA}_{\cdot,k}
}.
$$
Therefore,
$$
\E[s_\fG(z)]
=
\frac{1}{-z-\E[s_{\fG}(z)]}
+\delta_n,
$$
where
$$
\delta_n
=
\frac{1}{p}\sum_{k=1}^p
\E\Bigl[
\frac{\varepsilon_k}{
(z+\E[s_\fG(z)])(-z-\E[s_{\fG}(z)]+\varepsilon_k)
}
\Bigr],
$$
with
$$
\varepsilon_k
=
\E[s_{\fG}(z)]
-
y_n^{-1}M^{-2}\wt{\mA}_{\cdot,k}^\tp
\wt{\fA}_{-k}
\bigl(y_n^{-1/2}D_0(\fS_{-k})-z\bigr)^{-1}
\wt{\fA}_{-k}^\tp \wt{\mA}_{\cdot,k}.
$$
It therefore suffices to prove
\bea\label{Eq-limit-delta}
\lim_{n\to\infty}\delta_n=0.
\eea
Write $\delta_n=J_1+J_2$, where
$$
J_1
=
-\frac{1}{p}\sum_{k=1}^p
\frac{\E[\varepsilon_k]}{(z+\E[s_{\fG}(z)])^2},
$$
and
$$
J_2
=
\frac{1}{p}\sum_{k=1}^p
\E\Bigl[
\frac{\varepsilon_k^2}{
(z+\E[s_\fG(z)])^2(-z-\E[s_{\fG}(z)]+\varepsilon_k)
}
\Bigr].
$$
    

Since $\wt{\mA}_{\cdot,k}\perp \wt{\fA}_{-k}$ and
\[
\E[\wt{\mA}_{\cdot,k}\wt{\mA}_{\cdot,k}^\tp]=\fI_M,
\]
we have
\bea
|\E[\varepsilon_k]|
&=
\Bigl|
\frac{1}{p}\E\big[\tr(\fG-z)^{-1}\big]
-
\E\Bigl[
y_n^{-1}M^{-2}\tr\Bigl(
\wt{\fA}_{-k}\bigl(y_n^{-1/2}D_0(\fS_{-k})-z\bigr)^{-1}\wt{\fA}_{-k}^\tp
\Bigr)
\Bigr]
\Bigr|\\
&=
\Bigl|
\frac{1}{p}\E\big[\tr(\fG-z)^{-1}\big]
-
\frac{1}{p}\E\Bigl[
\tr\Bigl(
\frac{1}{M}\wt{\fA}_{-k}^\tp\wt{\fA}_{-k}
\bigl(y_n^{-1/2}D_0(\fS_{-k})-z\bigr)^{-1}
\Bigr)
\Bigr]
\Bigr|\\
&\leq
E_{k1}+E_{k2},
\eea
where
$$
E_{k1}
=
\frac{1}{p}
\Bigl|
\E\big[\tr(\fG-z)^{-1}\big]
-
\E\Big[\tr\bigl(y_n^{-1/2}D_0(\fS_{-k})-z\bigr)^{-1}\Big]
\Bigr|,
$$
and
$$
E_{k2}
=
\frac{1}{p}
\Bigl|
\E\Bigl[
\tr\Bigl(
(\fS_{-k}-\fI_{p-1})
\bigl(y_n^{-1/2}D_0(\fS_{-k})-z\bigr)^{-1}
\Bigr)
\Bigr]
\Bigr|.
$$
By \citet[equation~A.12]{bai2010spectral},
$$
E_{k1}\leq \frac{1}{pv}.
$$

To control $E_{k2}$, we need the following lemma.

\begin{lem}\label{lemma-Delta-Fnorm} $\E\big[\|\fS_{-k}-\fI_{p-1}\|_{\operatorname{F}}^2\big]\lesssim 1$.
\end{lem}

Then, by Hölder's inequality for Schatten norms,
\bea
E_{k2}
&\leq
\frac{1}{p}\E\Bigl[
\Bigl\|
(\fS_{-k}-\fI_{p-1})
\bigl(y_n^{-1/2}D_0(\fS_{-k})-z\bigr)^{-1}
\Bigr\|_{S_1}
\Bigr]\\
&\leq
\frac{1}{p}\E\Bigl[
\|\fS_{-k}-\fI_{p-1}\|_{S_1}
\cdot
\Bigl\|\bigl(y_n^{-1/2}D_0(\fS_{-k})-z\bigr)^{-1}\Bigr\|_{\mathrm{op}}
\Bigr]\\
&\leq
\frac{1}{v\sqrt{p}}\E\big[\|\fS_{-k}-\fI_{p-1}\|_{\mathrm{F}}\big]\\
&\leq
\frac{1}{v\sqrt{p}}
\sqrt{\E\Bigl[\|\fS_{-k}-\fI_{p-1}\|_{\mathrm{F}}^2\Bigr]}\\
&\lesssim
\frac{1}{\sqrt{p}}.
\eea
Hence,
$$
|\E[\varepsilon_k]|
\lesssim
\frac{1}{\sqrt{p}}
\to 0.
$$
Since
$$
|z+\E[s_\fG(z)]|
\geq
\Im\bigl(z+\E[s_\fG(z)]\bigr)
\geq v,
$$
we obtain
$$
J_1\to 0.
$$

For $J_2$, note that
\bea
&\big|-z-\E[s_\fG(z)]+\varepsilon_k\big|\\
=&
\Bigl|
-z-y_n^{-1}M^{-2}\wt{\mA}_{\cdot,k}^\tp
\wt{\fA}_{-k}
\bigl(y_n^{-1/2}D_0(\fS_{-k})-z\bigr)^{-1}
\wt{\fA}_{-k}^\tp \wt{\mA}_{\cdot,k}
\Bigr|\\
\geq&
\Im\Bigl(
z+y_n^{-1}M^{-2}\wt{\mA}_{\cdot,k}^\tp
\wt{\fA}_{-k}
\bigl(y_n^{-1/2}D_0(\fS_{-k})-z\bigr)^{-1}
\wt{\fA}_{-k}^\tp \wt{\mA}_{\cdot,k}
\Bigr)\\
\geq&
v
+
vy_n^{-1}M^{-2}\wt{\mA}_{\cdot,k}^\tp
\wt{\fA}_{-k}
\bigl(y_n^{-1/2}D_0(\fS_{-k})-z\bigr)^{-1}
\bigl(y_n^{-1/2}D_0(\fS_{-k})-\bar z\bigr)^{-1}
\wt{\fA}_{-k}^\tp \wt{\mA}_{\cdot,k}\\
\geq&
v.
\eea
    
Therefore,
$$
|J_2|
\leq
\frac{1}{pv^3}\sum_{k=1}^p \E[|\varepsilon_k|^2]
\lesssim
J_{21}+J_{22}+J_{23},
$$
where
\bea\nonumber
J_{21}
&=
\frac{1}{pv^3}\sum_{k=1}^p
\E\Big[|\varepsilon_k-\E[\varepsilon_k\mid\wt{\fA}_{-k}]|^2\Big],\\
J_{22}
&=
\frac{1}{pv^3}\sum_{k=1}^p
\E\Big[|\E[\varepsilon_k\mid\wt{\fA}_{-k}]-\E[\varepsilon_k]|^2\Big],\\
J_{23}
&=
\frac{1}{pv^3}\sum_{k=1}^p
|\E[\varepsilon_k]|^2.
\eea
For $J_{21}$, let
$$
\fB_k
:=
\wt{\fA}_{-k}
\bigl(y_n^{-1/2}D_0(\fS_{-k})-z\bigr)^{-1}
\wt{\fA}_{-k}^\tp
\in\R^{M\times M}.
$$
Then
$$
|\varepsilon_k-\E[\varepsilon_k\mid\wt{\fA}_{-k}]|
=
\frac{1}{pM}
\bigl|
\wt{\mA}_{\cdot,k}^\tp \fB_k \wt{\mA}_{\cdot,k}
-\tr(\fB_k)
\bigr|.
$$
Hence,
$$
|J_{21}|
\leq
\frac{1}{p^3M^2v^3}\sum_{k=1}^p
\E\Bigl[
\bigl|
\wt{\mA}_{\cdot,1}^\tp \fB_1 \wt{\mA}_{\cdot,1}
-\tr(\fB_1)
\bigr|^2
\Bigr].
$$

\begin{lem}\label{Lemma-Quadratic-form-control} We have
\bea\label{eq-quadratic-form-concentration}
\lim_{n\to\infty}
\frac{1}{p^3M^2}
\sum_{k=1}^p
\E\Bigl[
\bigl|
\wt{\mA}_{\cdot,1}^\tp \fB_1 \wt{\mA}_{\cdot,1}
-\tr(\fB_1)
\bigr|^2
\Bigr]
=
0.
\eea
\end{lem}

By Lemma~\ref{Lemma-Quadratic-form-control}, we have $J_{21}\to 0$.

For $J_{22}$, we first observe that
\bea\nonumber
\E\Big[|\E[\varepsilon_k\mid\wt{\fA}_{-k}]-\E[\varepsilon_k]|^2\Big]
&=
\E\Bigl[
\frac{1}{p^2M^2}\big|\tr\fB_k-\E[\tr\fB_k]\big|^2
\Bigr].
\eea
Note that
$$
\frac{1}{M}\tr\fB_k
=
\tr\bigl(y_n^{-1/2}D_0(\fS_{-k})-z\bigr)^{-1}
+
\tr\Bigl(
(\fS_{-k}-\fI_{p-1})
\bigl(y_n^{-1/2}D_0(\fS_{-k})-z\bigr)^{-1}
\Bigr),
$$
and, using $|a+b|^2\leq 2(a^2+b^2)$, we obtain
$$
J_{22}\leq J_{221}+J_{222},
$$
where
\bea\nonumber\footnotesize
J_{221}
&=
\frac{2}{p^3v^3}\sum_{k=1}^p
\E\Bigl[
\Bigl|
\tr\bigl(y_n^{-1/2}D_0(\fS_{-k})-z\bigr)^{-1}
-
\E\Big[\tr\bigl(y_n^{-1/2}D_0(\fS_{-k})-z\bigr)^{-1}\Big]
\Bigr|^2
\Bigr],\\
J_{222}
&=
\frac{2}{p^3v^3}\sum_{k=1}^p
\E\Bigl[
\Bigl|
\tr\Bigl(
(\fS_{-k}-\fI_{p-1})
\bigl(y_n^{-1/2}D_0(\fS_{-k})-z\bigr)^{-1}
\Bigr)
-
\E\Bigl[
\tr\Bigl(
(\fS_{-k}-\fI_{p-1})
\bigl(y_n^{-1/2}D_0(\fS_{-k})-z\bigr)^{-1}
\Bigr)
\Bigr]
\Bigr|^2
\Bigr].
\eea

For $J_{221}$, define $\rG_0:=\{\emptyset\}$ and, for $j\in[p]$ with $j\neq k$,
\[
\rG_j:=\sigma\Bigl(\{\mA_{\cdot,i}:1\leq i\leq j,\ i\neq k\}\Bigr),
\qquad
\E_{\rG_j}[\cdot]:=\E[\cdot\mid \rG_j].
\]
Then
$$
\tr\bigl(y_n^{-1/2}D_0(\fS_{-k})-z\bigr)^{-1}
-
\E\Big[\tr\bigl(y_n^{-1/2}D_0(\fS_{-k})-z\bigr)^{-1}\Big]
=
\sum_{j\in[p],\,j\neq k} W_{jk},
$$
where
\bea\nonumber
W_{jk}
:=
&(\E_{\rG_j}-\E_{\rG_{j-1}})
\Bigl[
\tr\bigl(y_n^{-1/2}D_0(\fS_{-k})-z\bigr)^{-1}
\Bigr]\\
=
&(\E_{\rG_j}-\E_{\rG_{j-1}})
\Bigl[
\tr\bigl(y_n^{-1/2}D_0(\fS_{-k})-z\bigr)^{-1}
-
\tr\bigl(y_n^{-1/2}D_0(\fS_{-kj})-z\bigr)^{-1}
\Bigr].
\eea
Here
\[
\fS_{-kj}:=\frac{1}{M}\wt{\fA}_{-kj}^\tp\wt{\fA}_{-kj},
\]
and $\wt{\fA}_{-kj}$ is the $M\times (p-2)$ matrix obtained by removing $\wt{\mA}_{\cdot,j}$ and $\wt{\mA}_{\cdot,k}$ from $\wt{\fA}$. By Lemma~\ref{lem-trace-resolvant-minor},
$$
\Bigl|
\tr\bigl(y_n^{-1/2}D_0(\fS_{-k})-z\bigr)^{-1}
-
\tr\bigl(y_n^{-1/2}D_0(\fS_{-kj})-z\bigr)^{-1}
\Bigr|
\leq
\frac{1}{v}.
$$
Hence $|W_{jk}|\leq 2v^{-1}$. By Lemma~\ref{lem-burkholder},
$$
\E\Bigl[
\Bigl|\sum_{j\in[p],\,j\neq k} W_{jk}\Bigr|^2
\Bigr]
\leq
C_2\E\Bigl[
\sum_{j\in[p],\,j\neq k} |W_{jk}|^2
\Bigr]
\leq
\frac{4C_2p}{v^2},
$$
where $C_2>0$ is a constant. Therefore,
\bea\nonumber
J_{221}
&\leq
\frac{2}{p^3v^3}\sum_{k=1}^p
\E\Bigl[
\Bigl|\sum_{j\in[p],\,j\neq k} W_{jk}\Bigr|^2
\Bigr]
\lesssim
\frac{1}{p}
\to 0.
\eea

 For $J_{222}$, using $|X-\E[X]|^2\leq 2|X|^2+2|\E[X]|^2$, we obtain
$$
\E\big[|X-\E[X]|^2\big]
\leq
2\E\big[|X|^2\big]+2|\E[X]|^2
\leq
4\E\big[|X|^2\big].
$$
Therefore, it suffices to control the second moment of
$$
\tr\Bigl(
(\fS_{-k}-\fI_{p-1})
\bigl(y_n^{-1/2}D_0(\fS_{-k})-z\bigr)^{-1}
\Bigr).
$$
By $|\tr(\fA)|\leq \|\fA\|_{S_1}$ and Hölder's inequality for Schatten norms,
\bea\nonumber
&\E\Bigl[
\Bigl|
\tr\Bigl(
(\fS_{-k}-\fI_{p-1})
\bigl(y_n^{-1/2}D_0(\fS_{-k})-z\bigr)^{-1}
\Bigr)
\Bigr|^2
\Bigr]\\
&\leq
\E\Bigl[
\|\fS_{-k}-\fI_{p-1}\|_{S_1}^2
\cdot
\Bigl\|\bigl(y_n^{-1/2}D_0(\fS_{-k})-z\bigr)^{-1}\Bigr\|_{\mathrm{op}}^2
\Bigr]\\
&\leq
\frac{1}{v^2}
\E\Bigl[
\|\fS_{-k}-\fI_{p-1}\|_{S_1}^2
\Bigr]\\
&\leq
\frac{p}{v^2}
\E\Bigl[
\|\fS_{-k}-\fI_{p-1}\|_{\mathrm{F}}^2
\Bigr]\\
&\lesssim
p,
\eea
where the last step follows from Lemma~\ref{lemma-Delta-Fnorm}. Hence
\bea\nonumber
J_{222}
&\lesssim
\frac{1}{p^3v^3}\sum_{k=1}^p
\E\Bigl[
\Bigl|
\tr\Bigl(
(\fS_{-k}-\fI_{p-1})
\bigl(y_n^{-1/2}D_0(\fS_{-k})-z\bigr)^{-1}
\Bigr)
\Bigr|^2
\Bigr]\\
&\lesssim
\frac{1}{p^3v^3}\cdot p\cdot p
\to 0.
\eea
Therefore, $J_{22}\to 0$.

For $J_{23}$, we have already shown that $|\E[\varepsilon_k]|\lesssim p^{-1/2}$. Hence
$$
J_{23}\lesssim p^{-1}\to 0.
$$

We conclude that $J_2\to 0$. Therefore, \eqref{Eq-limit-delta} holds, which completes the proof.

\subsection{Proof of the auxiliary lemmas}\label{section-proof-of-lemmas}

\subsubsection{Proof of Lemma \ref{lemma-cross-moments}}

By definition,
\bea\nonumber
\E\Bigl[\prod_{q=1}^Q \wt{A}_{\mi_q,1}\Bigr]
&=
\sigma_T^{-Q}
\E\Bigl[
\prod_{q=1}^Q
\Bigl(
\sum_{r=1}^T \lambda_r \psi_r(X_{i_1^q,1})\psi_r(X_{i_2^q,1})
\Bigr)
\Bigr]\\
&=
\sigma_T^{-Q}
\sum_{r_1,\ldots,r_Q=1}^T
\lambda_{r_1}\cdots\lambda_{r_Q}
\E\Bigl[
\psi_{r_1}(X_{i_1^1,1})\psi_{r_1}(X_{i_2^1,1})
\cdots
\psi_{r_Q}(X_{i_1^Q,1})\psi_{r_Q}(X_{i_2^Q,1})
\Bigr].
\eea
Without loss of generality, suppose that $i_1^1$ appears exactly once among the $2Q$ indices. Then, since
\[
\E[\psi_r(X_{i_1^1,1})]=0,
\]
it follows that for all $r_1,\ldots,r_Q\in[T]$,
\begin{align*}
&\E\Bigl[
\psi_{r_1}(X_{i_1^1,1})\psi_{r_1}(X_{i_2^1,1})
\cdots
\psi_{r_Q}(X_{i_1^Q,1})\psi_{r_Q}(X_{i_2^Q,1})
\Bigr]\\
&\qquad=
\E\big[\psi_{r_1}(X_{i_1^1,1})\big]\,
\E\Bigl[
\psi_{r_1}(X_{i_2^1,1})
\cdots
\psi_{r_Q}(X_{i_1^Q,1})\psi_{r_Q}(X_{i_2^Q,1})
\Bigr]
=0.
\end{align*}
Hence,
$$
\E\Bigl[\prod_{q=1}^Q \wt{A}_{\mi_q,1}\Bigr]=0.
$$

\subsubsection{Proof of Lemma \ref{lemma-Delta-Fnorm}}

Since $\wt{A}_{\mu,1}\perp \wt{A}_{\nu,1}$ whenever $\mu\cap\nu=\emptyset$, we have
\bea\nonumber
\E\Bigl[\Bigl(\frac{1}{M}\wt{\mA}_{\cdot,1}^\tp\wt{\mA}_{\cdot,1}-1\Bigr)^2\Bigr]
&=
\frac{1}{M^2}
\E\Bigl[
\Bigl(\sum_{\mu\in\cM_2}(\wt{A}_{\mu,1}^2-1)\Bigr)^2
\Bigr]\\
&=
\frac{1}{M^2}
\sum_{\mu,\nu\in\cM_2}
\E\Bigl[(\wt{A}_{\mu,1}^2-1)(\wt{A}_{\nu,1}^2-1)\Bigr]
\cdot
\ind(\mu\cap\nu\neq\emptyset)\\
&\lesssim
\frac{1}{M^2}\cdot n^3\\
&\lesssim
\frac{1}{n}.
\eea
Moreover, since $\wt{\mA}_{\cdot,1}\perp \wt{\mA}_{\cdot,2}$,
\bea\nonumber
\E\Bigl[\Bigl(\frac{1}{M}\wt{\mA}_{\cdot,1}^\tp\wt{\mA}_{\cdot,2}\Bigr)^2\Bigr]
&=
\frac{1}{M^2}
\sum_{\mu,\nu\in\cM_2}
\E[\wt{A}_{\mu,1}\wt{A}_{\mu,2}\wt{A}_{\nu,1}\wt{A}_{\nu,2}]\\
&=
\frac{1}{M^2}
\sum_{\mu,\nu\in\cM_2}
\E[\wt{A}_{\mu,1}\wt{A}_{\nu,1}]
\E[\wt{A}_{\mu,2}\wt{A}_{\nu,2}]\\
&=
\frac{1}{M}.
\eea
Hence,
\bea\nonumber
\E[\|\fS-\fI_p\|_{\mathrm{F}}^2]
&=
p\,\E\Bigl[\Bigl(\frac{1}{M}\wt{\mA}_{\cdot,1}^\tp\wt{\mA}_{\cdot,1}-1\Bigr)^2\Bigr]
+
p(p-1)\,\E\Bigl[\Bigl(\frac{1}{M}\wt{\mA}_{\cdot,1}^\tp\wt{\mA}_{\cdot,2}\Bigr)^2\Bigr]
\lesssim
1.
\eea
Similarly, we also have
\[
\E[\|\fS_{-k}-\fI_{p-1}\|_{\mathrm{F}}^2]\lesssim 1.
\]

\subsubsection{Proof of Lemma \ref{Lemma-Quadratic-form-control}}

Define
$$
\fQ_k:=\bigl(y_n^{-1/2}D_0(\fS_{-k})-z\bigr)^{-1},
\qquad
\fD_k:=\fS_{-k}-\fI_{p-1}.
$$
Then $\fB_k=\wt{\fA}_{-k}\fQ_k\wt{\fA}_{-k}^\tp$, and
\bea\nonumber
\frac{1}{M^2}\E\big[\|\fB_k\|_{\mathrm{F}}^2\big]
&=
\frac{1}{M^2}\E\Bigl[
\tr\Bigl(
\wt{\fA}_{-k}\fQ_k\wt{\fA}_{-k}^\tp
\wt{\fA}_{-k}\fQ_k^*\wt{\fA}_{-k}^\tp
\Bigr)
\Bigr]\\
&=
\E\Bigl[
\tr\Bigl(
(\fI+\fD_k)\fQ_k(\fI+\fD_k)\fQ_k^*
\Bigr)
\Bigr]\\
&=
\E\big[\tr(\fQ_k\fQ_k^*)\big]
+\E\big[\tr(\fD_k\fQ_k\fQ_k^*)\big]
+\E\big[\tr(\fQ_k\fD_k\fQ_k^*)\big]
+\E\big[\tr(\fD_k\fQ_k\fD_k\fQ_k^*)\big].
\eea
Using $\E[\|\fD_k\|_{\mathrm{F}}^2]\lesssim 1$ from Lemma~\ref{lemma-Delta-Fnorm}, the bound $\|\fQ_k\|_{\mathrm{op}}\leq v^{-1}$, the inequality $|\tr(\fA)|\leq \|\fA\|_{S_1}$, and Hölder's inequality for Schatten norms, we obtain
$$
\E[\tr(\fQ_k\fQ_k^*)]\lesssim p,
$$
$$
|\E[\tr(\fD_k\fQ_k\fQ_k^*)]|
\leq
\E[\|\fD_k\fQ_k\fQ_k^*\|_{S_1}]
\lesssim
\E[\|\fD_k\|_{S_1}\|\fQ_k\fQ_k^*\|_{S_\infty}]
\lesssim
\sqrt{p}\,\E[\|\fD_k\|_{\mathrm{F}}]
\lesssim
\sqrt{p},
$$
$$
|\E[\tr(\fQ_k\fD_k\fQ_k^*)]|
=
|\E[\tr(\fD_k\fQ_k^*\fQ_k)]|
\lesssim
\E[\|\fD_k\|_{S_1}]
\lesssim
\sqrt{p},
$$
and
\bea\nonumber
|\E[\tr(\fD_k\fQ_k\fD_k\fQ_k^*)]|
&\leq
\E[\|\fD_k\fQ_k\fD_k\fQ_k^*\|_{S_1}]\\
&\leq
\E[\|\fD_k\|_{S_2}\|\fQ_k\|_{S_\infty}\|\fD_k\|_{S_2}\|\fQ_k\|_{S_\infty}]\\
&\lesssim
\E[\|\fD_k\|_{\mathrm{F}}^2]\\
&\lesssim
1.
\eea
Hence,
\bea\nonumber
\E[\|\fB_k\|_{\mathrm{F}}^2]\lesssim pM^2.
\eea
Write $\fB_k=(b_{\mu\nu})_{\mu,\nu\in\cM_2}$. Since
$$
\E[\|\fB_k\|_{\mathrm{F}}^2]
=
M\,\E[|b_{\mu\mu}|^2]
+
M(M-1)\,\E[|b_{\mu\nu}|^2\ind(\mu\neq \nu)],
$$
we have
\bea\label{eq-b_munu}
\E[|b_{\mu\nu}|^2\ind(\mu\neq \nu)]\lesssim p,
\eea
and, by the Cauchy--Schwarz inequality,
\bea\label{eq-b_munu-b_nukappa}
\E[|b_{\mu\nu}b_{\kappa\tau}|\ind(\mu\neq \nu,\kappa\neq \tau)]\lesssim p.
\eea
Let $\wt{\mA}^{(-k)}_{\mu,\cdot}$ denote the $\mu$th row of $\wt{\fA}_{-k}$. Then
\bea\label{eq-bmumu}
\E[|b_{\mu\mu}|^2]
&=
\E\Bigl[
\wt{\mA}^{(-k)}_{\mu,\cdot}\fQ_k\wt{\mA}^{(-k)\tp}_{\mu,\cdot}
\wt{\mA}^{(-k)}_{\mu,\cdot}\fQ_k^*\wt{\mA}^{(-k)\tp}_{\mu,\cdot}
\Bigr]
\\
&\leq
\E\Bigl[
\Bigl|
\tr\Bigl(
\wt{\mA}^{(-k)\tp}_{\mu,\cdot}\wt{\mA}^{(-k)}_{\mu,\cdot}
\fQ_k
\wt{\mA}^{(-k)\tp}_{\mu,\cdot}\wt{\mA}^{(-k)}_{\mu,\cdot}
\fQ_k^*
\Bigr)
\Bigr|
\Bigr]\\
&\leq
\E\Bigl[
\Bigl\|
\wt{\mA}^{(-k)\tp}_{\mu,\cdot}\wt{\mA}^{(-k)}_{\mu,\cdot}
\fQ_k
\wt{\mA}^{(-k)\tp}_{\mu,\cdot}\wt{\mA}^{(-k)}_{\mu,\cdot}
\fQ_k^*
\Bigr\|_{S_1}
\Bigr]\\
&\leq
\frac{1}{v^2}
\E\Bigl[
\Bigl\|
\wt{\mA}^{(-k)\tp}_{\mu,\cdot}\wt{\mA}^{(-k)}_{\mu,\cdot}
\Bigr\|_{S_2}^2
\Bigr]
\\
&=
\frac{1}{v^2}
\E\Bigl[
\Bigl(
\sum_{\ell=1}^{p-1} (\wt{\mA}^{(-k)}_{\mu,\ell})^2
\Bigr)^2
\Bigr]\\
&\lesssim
p^2.
\eea

We are now ready to prove \eqref{eq-quadratic-form-concentration}. We have
$$
\E\Bigl[\Bigl(\wt{\mA}_{\cdot,k}^\tp\fB_{k}\wt{\mA}_{\cdot,k}-\tr\fB_k\Bigr)^2\Bigr]
\leq
2\E\Bigl[\Bigl(\sum_{\mu\neq\nu}b_{\mu\nu}\wt{A}_{\mu,k}\wt{A}_{\nu,k}\Bigr)^2\Bigr]
+
2\E\Bigl[\Bigl(\sum_{\mu\in\cM_2}b_{\mu\mu}(\wt{A}_{\mu,k}^2-1)\Bigr)^2\Bigr].
$$
Since $\wt{\mA}_{\cdot,k}\perp \fB_k$, Lemma~\ref{lemma-cross-moments} and the uniform boundedness of $(\wt{A}_{\mu,k})$ imply that
\bea\nonumber
\E\Bigl[\Bigl(\sum_{\mu\neq\nu}b_{\mu\nu}\wt{A}_{\mu,k}\wt{A}_{\nu,k}\Bigr)^2\Bigr]
&\leq
\sum_{\mu\neq\nu,\ \kappa\neq\tau}
\E[|b_{\mu\nu}b_{\kappa\tau}|]\cdot
\Bigl|
\E\big[\wt{A}_{\mu,k}\wt{A}_{\nu,k}\wt{A}_{\kappa,k}\wt{A}_{\tau,k}\big]
\Bigr|\\
&\lesssim
F_1+F_2+F_3,
\eea
where
\bea\nonumber
F_1
&=
\sum_{\mu,\nu,\kappa,\tau}
\E\Bigl[
|b_{\mu\nu}|^2\cdot \ind(\mu=\kappa,\ \nu=\tau,\ \mu\neq\nu)
\Bigr],\\
F_2
&=
\sum_{\mu,\nu,\kappa,\tau}
\E\Bigl[
|b_{\mu\nu}b_{\kappa\tau}|\cdot
\ind(\tau=\nu,\ \mu\cap\kappa\neq\emptyset,\ \mu\cap\nu\neq\emptyset,\ \nu\cap\tau\neq\emptyset,\ \mu\neq\nu\neq\kappa)
\Bigr],\\
F_3
&=
\sum_{\mu,\nu,\kappa,\tau}
\E\Bigl[
|b_{\mu\nu}b_{\kappa\tau}|\cdot
\ind(\mu\cap\nu\neq\emptyset,\ \mu\cap\kappa\neq\emptyset,\ \nu\cap\tau\neq\emptyset,\ \tau\cap\kappa\neq\emptyset,\ \mu\neq\nu\neq\kappa\neq\tau)
\Bigr].
\eea
By \eqref{eq-b_munu} and \eqref{eq-b_munu-b_nukappa},
$$
F_1
=
\sum_{\mu,\nu}
\E[|b_{\mu\nu}|^2\ind(\mu\neq\nu)]
\lesssim
pM^2,
$$
$$
F_2
\lesssim
p\cdot
\sum_{\mu\neq\nu\neq\kappa}
\ind(\mu\cap\kappa\neq\emptyset,\ \mu\cap\nu\neq\emptyset,\ \nu\cap\kappa\neq\emptyset)
\lesssim
p^4,
$$
$$
F_3
\lesssim
p\cdot
\sum_{\mu\neq\nu\neq\kappa\neq\tau}
\ind(\mu\cap\nu\neq\emptyset,\ \mu\cap\kappa\neq\emptyset,\ \nu\cap\tau\neq\emptyset,\ \tau\cap\kappa\neq\emptyset)
\lesssim
p^5.
$$
Hence,
\bea\label{eq-quadratic-form-control-dediagonal}
\E\Bigl[\Bigl(\sum_{\mu\neq\nu}b_{\mu\nu}\wt{A}_{\mu,k}\wt{A}_{\nu,k}\Bigr)^2\Bigr]
\lesssim
p^5.
\eea

Similarly, using the fact that $\wt{A}_{\mu,k}\perp \wt{A}_{\nu,k}$ whenever $\mu\cap\nu=\emptyset$, we have
\bea
\E\Bigl[\Bigl(\sum_{\mu}b_{\mu\mu}(\wt{A}_{\mu,k}^2-1)\Bigr)^2\Bigr]
&\leq
\sum_{\mu,\nu}
\E[|b_{\mu\mu}b_{\nu\nu}|]\cdot
\Bigl|
\E\big[(\wt{A}_{\mu,k}^2-1)(\wt{A}_{\nu,k}^2-1)\big]
\Bigr|\\
&=
\sum_{\mu,\nu}
\E[|b_{\mu\mu}b_{\nu\nu}|]\cdot
\Bigl|
\E\big[(\wt{A}_{\mu,k}^2-1)(\wt{A}_{\nu,k}^2-1)\big]
\Bigr|\cdot
\ind(\mu\cap\nu\neq\emptyset)\\
&\lesssim
G_1+G_2,
\eea
where
$$
G_1=\sum_{\mu}\E[|b_{\mu\mu}|^2]~~~{\rm and}~~~
G_2=
\sum_{\mu,\nu}
\E\Bigl[
|b_{\mu\mu}b_{\nu\nu}|\ind(\mu\cap\nu\neq\emptyset,\ \mu\neq\nu)
\Bigr].
$$
Then, by \eqref{eq-bmumu} and the Cauchy--Schwarz inequality,
$$
G_1\lesssim p^2M~~~{\rm and}~~~
G_2
\lesssim
p^2\sum_{\mu\neq\nu}\ind(\mu\cap\nu\neq\emptyset)
\lesssim
p^5.
$$
Hence,
\bea\label{eq-quadratic-form-control-diagonal}
\E\Bigl[\Bigl(\sum_{\mu\in\cM_2}b_{\mu\mu}(\wt{A}_{\mu,k}^2-1)\Bigr)^2\Bigr]
\lesssim
p^5.
\eea
Combining \eqref{eq-quadratic-form-control-dediagonal} and \eqref{eq-quadratic-form-control-diagonal}, we conclude that
$$
\E\Bigl[\Bigl(\wt{\mA}_{\cdot,k}^\tp\fB_{k}\wt{\mA}_{\cdot,k}-\tr\fB_k\Bigr)^2\Bigr]
\lesssim
p^5,
$$
and hence \eqref{eq-quadratic-form-concentration} follows.

\section{Appendix}

\begin{lem}\label{lem-perturbation-frob-norm}\citep[Lemma~F.1]{lu2025equivalence}
Let $\fA,\fB$ be two $p\times p$ Hermitian matrices, and let $z\in\bC^+$ with $v=\Im z>0$. Then
$$
|s_{\fA}(z)-s_{\fB}(z)|
\leq
\frac{\|\fA-\fB\|_{\operatorname{F}}}{\sqrt{p}\,v^2}.
$$
\end{lem}

 \begin{lem}\label{lem-trace-resolvant-minor}\citep[Theorem~A.6]{bai2010spectral}
Let $\fA$ be a $p\times p$ Hermitian matrix, and let $\fA_k$ be the $(p-1)\times (p-1)$ matrix obtained by removing the $k$th row and column of $\fA$. For $z\in\bC^+$ with $v=\Im z>0$, we have
$$
\Bigl|
\tr\bigl(\fA-z\fI_p\bigr)^{-1}
-
\tr\bigl(\fA_k-z\fI_{p-1}\bigr)^{-1}
\Bigr|
\leq
\frac{1}{v}.
$$
\end{lem}

\begin{lem}\label{lem-burkholder}\citep[Lemma~2.12]{bai2010spectral}
Let $\{X_k\}_{k=1}^n$ be a complex martingale difference sequence with respect to the filtration $\{\rF_k\}_{k=0}^n$. Then, for $p>1$,
$$
\E\Bigl[\Bigl|\sum_{k=1}^n X_k\Bigr|^p\Bigr]
\leq
K_p\E\Bigl[\Bigl(\sum_{k=1}^n |X_k|^2\Bigr)^{p/2}\Bigr],
$$
where $K_p>0$ is a constant.
\end{lem}
    
{
\bibliographystyle{apalike}
\bibliography{reference}
}

\end{document}